# COMMON FUNCTIONAL PRINCIPAL COMPONENTS[1]

By Michal Benko, Wolfgang Härdle and Alois Kneip

*Humboldt-Universität, Humboldt-Universität and Bonn Universität*

Functional principal component analysis (FPCA) based on the Karhunen–Loève decomposition has been successfully applied in many applications, mainly for one sample problems. In this paper we consider common functional principal components for two sample problems. Our research is motivated not only by the theoretical challenge of this data situation, but also by the actual question of dynamics of implied volatility (IV) functions. For different maturities the log-returns of IVs are samples of (smooth) random functions and the methods proposed here study the similarities of their stochastic behavior. First we present a new method for estimation of functional principal components from discrete noisy data. Next we present the two sample inference for FPCA and develop the two sample theory. We propose bootstrap tests for testing the equality of eigenvalues, eigenfunctions, and mean functions of two functional samples, illustrate the test-properties by simulation study and apply the method to the IV analysis.

**1. Introduction.** In many applications in biometrics, chemometrics, econometrics, etc., the data come from the observation of continuous phenomenons of time or space and can be assumed to represent a sample of i.i.d. smooth random functions $X_1(t), \ldots, X_n(t) \in L^2[0, 1]$. Functional data analysis has received considerable attention in the statistical literature during the last decade. In this context functional principal component analysis (FPCA) has proved to be a key technique. An early reference is Rao (1958), and important methodological contributions have been given by various authors. Case studies and references, as well as methodological and algorithmical details, can be found in the books by Ramsay and Silverman (2002, 2005) or Ferraty and Vieu (2006).

Received January 2006; revised February 2007.

[1]Supported by the Deutsche Forschungsgemeinschaft and the Sonderforschungsbereich 649 "Ökonomisches Risiko."

The well-known Karhunen–Loève (KL) expansion provides a basic tool to describe the distribution of the random functions $X_i$ and can be seen as the theoretical basis of FPCA. For $v, w \in L^2[0,1]$, let $\langle v, w \rangle = \int_0^1 v(t)w(t)\,dt$, and let $\| \cdot \| = \langle \cdot, \cdot \rangle^{1/2}$ denote the usual $L^2$-norm. With $\lambda_1 \geq \lambda_2 \geq \cdots$ and $\gamma_1, \gamma_2, \ldots$ denoting eigenvalues and corresponding orthonormal eigenfunctions of the covariance operator $\Gamma$ of $X_i$, we obtain $X_i = \mu + \sum_{r=1}^{\infty} \beta_{ri}\gamma_r, i = 1, \ldots, n$, where $\mu = \mathrm{E}(X_i)$ is the mean function and $\beta_{ri} = \langle X_i - \mu, \gamma_r \rangle$ are (scalar) factor loadings with $\mathrm{E}(\beta_{ri}^2) = \lambda_r$. Structure and dynamics of the random functions can be assessed by analyzing the "functional principal components" $\gamma_r$, as well as the distribution of the factor loadings. For a given functional sample, the unknown characteristics $\lambda_r, \gamma_r$ are estimated by the eigenvalues and eigenfunctions of the empirical covariance operator $\hat{\Gamma}_n$ of $X_1, \ldots, X_n$. Note that an eigenfunction $\gamma_r$ is identified (up to sign) only if the corresponding eigenvalue $\lambda_r$ has multiplicity one. This therefore establishes a necessary regularity condition for any inference based on an estimated functional principal component $\hat{\gamma}_r$ in FPCA. Signs are arbitrary ($\gamma_r$ and $\beta_{ri}$ can be replaced by $-\gamma_r$ and $-\beta_{ri}$) and may be fixed by a suitable standardization. More detailed discussion on this topic and precise assumptions can be found in Section 2.

In many important applications a small number of functional principal components will suffice to approximate the functions $X_i$ with a high degree of accuracy. Indeed, FPCA plays a much more central role in functional data analysis than its well-known analogue in multivariate analysis. There are two major reasons. First, distributions on function spaces are complex objects, and the Karhunen–Loève expansion seems to be the only practically feasible way to access their structure. Second, in multivariate analysis a substantial interpretation of principal components is often difficult and has to be based on vague arguments concerning the correlation of principal components with original variables. Such a problem does not at all exists in the functional context, where $\gamma_1(t), \gamma_2(t), \ldots$ are functions representing the major modes of variation of $X_i(t)$ over $t$.

In this paper we consider inference and tests of hypotheses on the structure of functional principal components. Motivated by an application to implied volatility analysis, we will concentrate on the two sample case. A central point is the use of bootstrap procedures. We will show that the bootstrap methodology can also be applied to functional data.

In Section 2 we start by discussing one-sample inference for FPCA. Basic results on asymptotic distributions have already been derived by Dauxois, Pousse and Romain (1982) in situations where the functions are directly observable. Hall and Hosseini-Nasab (2006) develop asymptotic Taylor expansions of estimated eigenfunctions in terms of the difference $\hat{\Gamma}_n - \Gamma$.

Without deriving rigorous theoretical results, they also provide some qualitative arguments as well as simulation results motivating the use of bootstrap in order to construct confidence regions for principal components.

In practice, the functions of interest are often not directly observed, but are regression curves which have to be reconstructed from discrete, noisy data. In this context the standard approach is to first estimate individual functions nonparametrically (e.g., by B-splines) and then to determine principal components of the resulting estimated empirical covariance operator—see Besse and Ramsay (1986), Ramsay and Dalzell (1991), among others. Approaches incorporating a smoothing step into the eigenanalysis have been proposed by Rice and Silverman (1991), Pezzulli and Silverman (1993) or Silverman (1996). Robust estimation of principal components has been considered by Lacontore et al. (1999). Yao, Müller and Wang (2005) and Hall, Müller and Wang (2006) propose techniques based on nonparametric estimation of the covariance function $\mathrm{E}[\{X_i(t) - \mu(t)\}\{X_i(s) - \mu(s)\}]$ which can also be applied if there are only a few scattered observations per curve.

Section 2.1 presents a new method for estimation of functional principal components. It consists in an adaptation of a technique introduced by Kneip and Utikal (2001) for the case of density functions. The key-idea is to represent the components of the Karhunen–Loève expansion in terms of an $(L^2)$ scalar-product matrix of the sample. We investigate the asymptotic properties of the proposed method. It is shown that under mild conditions the additional error caused by estimation from discrete, noisy data is first-order asymptotically negligible, and inference may proceed "as if" the functions were directly observed. Generalizing the results of Dauxois, Pousse and Romain (1982), we then present a theorem on the asymptotic distributions of the empirical eigenvalues and eigenfunctions. The structure of the asymptotic expansion derived in the theorem provides a basis to show consistency of bootstrap procedures.

Section 3 deals with two-sample inference. We consider two independent samples of functions $\{X_i^{(1)}\}_{i=1}^{n_1}$ and $\{X_i^{(2)}\}_{i=1}^{n_2}$. The problem of interest is to test in how far the distributions of these random functions coincide. The structure of the different distributions in function space can be accessed by means of the respective Karhunen–Loève expansions

$$X_i^{(p)} = \mu^{(p)} + \sum_{r=1}^{\infty} \beta_{ri}^{(p)} \gamma_r^{(p)}, \qquad p = 1, 2.$$

Differences in the distribution of these random functions will correspond to differences in the components of the respective KL expansions above. Without restriction, one may require that signs are such that $\langle \gamma_r^{(1)}, \gamma_r^{(2)} \rangle \geq 0$. Two sample inference for FPCA in general has not been considered in the literature so far. In Section 3 we define bootstrap procedures for testing the equality of mean functions, eigenvalues, eigenfunctions and eigenspaces. Consistency of the bootstrap is derived in Section 3.1, while Section 3.2 contains a simulation study providing insight into the finite sample performance of our tests.

It is of particular interest to compare the functional components characterizing the two samples. If these factors are "common," this means $\gamma_r :=$ $\gamma_r^{(1)} = \gamma_r^{(2)}$, then only the factor loadings $\beta_{ri}^{(p)}$ may vary across samples. This situation may be seen as a functional generalization of the concept of "common principal components" as introduced by Flury (1988) in multivariate analysis. A weaker hypothesis may only require equality of the eigenspaces spanned by the first $L \in \mathbb{N}$ functional principal components. [$\mathbb{N}$ denotes the set of all natural numbers $1, 2, \ldots$ $(0 \notin \mathbb{N})$]. If for both samples the common $L$-dimensional eigenspaces suffice to approximate the functions with high accuracy, then the distributions in function space are well represented by a low-dimensional factor model, and subsequent analysis may rely on comparing the multivariate distributions of the random vectors $(\beta_{r1}^{(p)}, \ldots, \beta_{rL}^{(p)})^\top$.

The idea of "common functional principal components" is of considerable importance in implied volatility (IV) dynamics. This application is discussed in detail in Section 4. Implied volatility is obtained from the pricing model proposed by Black and Scholes (1973) and is a key parameter for quoting options prices. Our aim is to construct low-dimensional factor models for the log-returns of the IV functions of options with different maturities. In our application the first group of functional observations—$\{X_i^{(1)}\}_{i=1}^{n_1}$, are log-returns on the maturity "1 month" (1M group) and second group—$\{X_i^{(2)}\}_{i=1}^{n_2}$, are log-returns on the maturity "3 months" (3M group).

The first three eigenfunctions (ordered with respect to the corresponding eigenvalues), estimated by the method described in Section 2.1, are plotted in Figure 1. The estimated eigenfunctions for both groups are of similar structure, which motivates a common FPCA approach. Based on discretized vectors of functional values, a (multivariate) common principal components analysis of implied volatilities has already been considered by Fengler, Härdle and Villa (2003). They rely on the methodology introduced by Flury (1988) which is based on maximum likelihood estimation under the assumption of multivariate normality. Our analysis overcomes the limitations of this approach by providing specific hypothesis tests in a fully functional setup. It will be shown in Section 4 that for both groups $L = 3$ components suffice to explain 98.2% of the variability of the sample functions. An application of the tests developed in Section 3 does not reject the equality of the corresponding eigenspaces.

## 2. Functional principal components and one sample inference.
In this section we will focus on one sample of i.i.d. smooth random functions $X_1, \ldots,$

$X_n \in L^2[0,1]$. We will assume a well-defined mean function $\mu = \mathrm{E}(X_i)$, as well as the existence of a continuous covariance function $\sigma(t,s) = \mathrm{E}[\{X_i(t) - \mu(t)\}\{X_i(s) - \mu(s)\}]$. Then $\mathrm{E}(\|X_i - \mu\|^2) = \int \sigma(t,t)\,dt < \infty$, and the covariance operator $\Gamma$ of $X_i$ is given by

$$(\Gamma v)(t) = \int \sigma(t,s)v(s)\,ds, \qquad v \in L^2[0,1].$$

The Karhunen–Loève decomposition provides a basic tool to describe the distribution of the random functions $X_i$. With $\lambda_1 \geq \lambda_2 \geq \cdots$ and $\gamma_1, \gamma_2, \ldots$ denoting eigenvalues and a corresponding complete orthonormal basis of eigenfunctions of $\Gamma$, we obtain

$$(1) \qquad X_i = \mu + \sum_{r=1}^{\infty} \beta_{ri}\gamma_r, \qquad i = 1, \ldots, n,$$

where $\beta_{ri} = \langle X_i - \mu, \gamma_r \rangle$ are uncorrelated (scalar) factor loadings with $\mathrm{E}(\beta_{ri}) = 0$, $\mathrm{E}(\beta_{ri}^2) = \lambda_r$ and $\mathrm{E}(\beta_{ri}\beta_{ki}) = 0$ for $r \neq k$. Structure and dynamics of the random functions can be assessed by analyzing the "functional principal components" $\gamma_r$, as well as the distribution of the factor loadings.

A discussion of basic properties of (1) can, for example, be found in Gihman and Skorohod (1973). Under our assumptions, the infinite sums in (1) converge with probability 1, and $\sum_{r=1}^{\infty} \lambda_r = \mathrm{E}(\|X_i - \mu\|^2) < \infty$. Smoothness of $X_i$ carries over to a corresponding degree of smoothness of $\sigma(t,s)$ and $\gamma_r$. If, with probability 1, $X_i(t)$ is twice continuously differentiable, then $\sigma$ as well as $\gamma_r$ are also twice continuously differentiable. The particular case of a Gaussian random function $X_i$ implies that the $\beta_{ri}$ are independent $N(0, \lambda_r)$-distributed random variables.

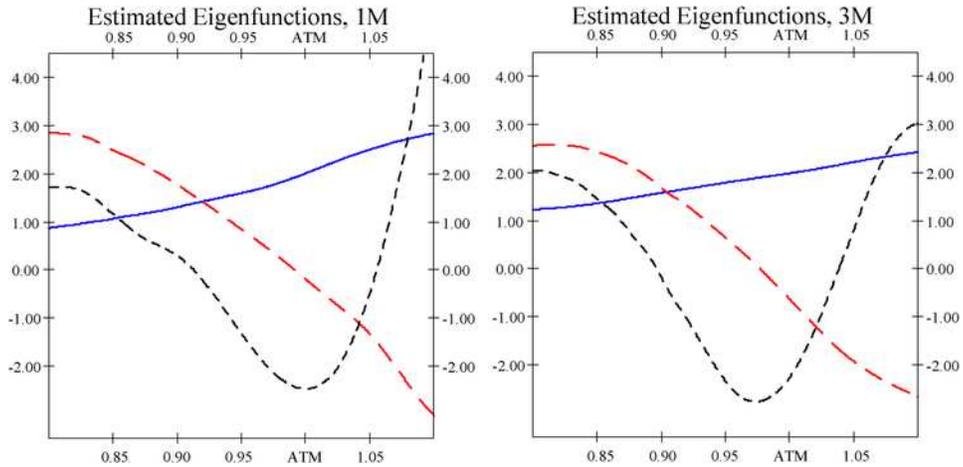

FIG. 1.   *Estimated eigenfunctions for 1M group in the left plot and 3M group in the right plot: solid—first function, dashed—second function, finely dashed—third function.*

An important property of (1) consists in the known fact that the first $L$ principal components provide a "best basis" for approximating the sample functions in terms of the integrated square error; see Ramsay and Silverman (2005), Section 6.2.3, among others. For any choice of $L$ orthonormal basis functions $v_1, \ldots, v_L$, the mean integrated square error

$$(2) \qquad \rho(v_1, \ldots, v_L) = \mathrm{E}\left( \left\| X_i - \mu - \sum_{r=1}^{L} \langle X_i - \mu, v_r \rangle v_r \right\|^2 \right)$$

is minimized by $v_r = \gamma_r$.

2.1. *Estimation of functional principal components.* For a given sample an empirical analog of (1) can be constructed by using eigenvalues $\hat{\lambda}_1 \geq \hat{\lambda}_2 \geq \cdots$ and orthonormal eigenfunctions $\hat{\gamma}_1, \hat{\gamma}_2, \ldots$ of the empirical covariance operator $\hat{\Gamma}_n$, where

$$(\hat{\Gamma}_n v)(t) = \int \hat{\sigma}(t, s) v(s) \, ds,$$

with $\bar{X} = n^{-1} \sum_{i=1}^n X_i$ and $\hat{\sigma}(t, s) = n^{-1} \sum_{i=1}^n \{X_i(t) - \bar{X}(t)\}\{X_i(s) - \bar{X}(s)\}$ denoting sample mean and covariance function. Then

$$(3) \qquad X_i = \bar{X} + \sum_{r=1}^{n} \hat{\beta}_{ri} \hat{\gamma}_r, \qquad i = 1, \ldots, n,$$

where $\hat{\beta}_{ri} = \langle \hat{\gamma}_r, X_i - \bar{X} \rangle$. We necessarily obtain $n^{-1} \sum_i \hat{\beta}_{ri} = 0$, $n^{-1} \sum_i \hat{\beta}_{ri} \hat{\beta}_{si} = 0$ for $r \neq s$, and $n^{-1} \sum_i \hat{\beta}_{ri}^2 = \hat{\lambda}_r$.

Analysis will have to concentrate on the leading principal components explaining the major part of the variance. In the following we will assume that $\lambda_1 > \lambda_2 > \cdots > \lambda_{r_0} > \lambda_{r_0+1}$, where $r_0$ denotes the maximal number of components to be considered. For all $r = 1, \ldots, r_0$, the corresponding eigenfunction $\gamma_r$ is then uniquely defined up to sign. Signs are arbitrary, decompositions (1) or (3) may just as well be written in terms of $-\gamma_r, -\beta_{ri}$ or $-\hat{\gamma}_r, -\hat{\beta}_{ri}$, and any suitable standardization may be applied by the statistician. In order to ensure that $\hat{\gamma}_r$ may be viewed as an estimator of $\gamma_r$ rather than of $-\gamma_r$, we will in the following only assume that signs are such that $\langle \gamma_r, \hat{\gamma}_r \rangle \geq 0$. More generally, any subsequent statement concerning differences of two eigenfunctions will be based on the condition of a nonnegative inner product. This does not impose any restriction and will go without saying.

The results of Dauxois, Pousse and Romain (1982) imply that, under regularity conditions, $\|\hat{\gamma}_r - \gamma_r\| = \mathcal{O}_p(n^{-1/2})$, $|\hat{\lambda}_r - \lambda_r| = \mathcal{O}_p(n^{-1/2})$, as well as $|\hat{\beta}_{ri} - \beta_{ri}| = \mathcal{O}_p(n^{-1/2})$ for all $r \leq r_0$.

However, in practice, the sample functions $X_i$ are often not directly observed, but have to be reconstructed from noisy observations $Y_{ij}$ at discrete

design points $t_{ik}$:

$$(4) \qquad Y_{ik} = X_i(t_{ik}) + \varepsilon_{ik}, \qquad k = 1, \ldots, T_i,$$

where $\varepsilon_{ik}$ are independent noise terms with $\mathrm{E}(\varepsilon_{ik}) = 0$, $\mathrm{Var}(\varepsilon_{ik}) = \sigma_i^2$.

Our approach for estimating principal components is motivated by the well-known duality relation between row and column spaces of a data matrix; see Härdle and Simar (2003), Chapter 8, among others. In a first step this approach relies on estimating the elements of the matrix:

$$(5) \qquad M_{lk} = \langle X_l - \bar{X}, X_k - \bar{X} \rangle, \qquad l, k = 1, \ldots, n.$$

Some simple linear algebra shows that all nonzero eigenvalues $\hat{\lambda}_1 \geq \hat{\lambda}_2 \cdots$ of $\hat{\Gamma}_n$ and $l_1 \geq l_2 \cdots$ of $M$ are related by $\hat{\lambda}_r = l_r/n$, $r = 1, 2, \ldots$. When using the corresponding orthonormal eigenvectors $p_1, p_2, \ldots$ of $M$, the empirical scores $\hat{\beta}_{ri}$, as well as the empirical eigenfunctions $\hat{\gamma}_r$, are obtained by $\hat{\beta}_{ri} = \sqrt{l_r} p_{ir}$ and

$$(6) \qquad \hat{\gamma}_r = \frac{1}{\sqrt{l_r}} \sum_{i=1}^{n} p_{ir}(X_i - \bar{X}) = \frac{1}{\sqrt{l_r}} \sum_{i=1}^{n} p_{ir} X_i.$$

The elements of $M$ are functionals which can be estimated with asympotically negligible bias and a parametric rate of convergence $T_i^{-1/2}$. If the data in (4) is generated from a balanced, equidistant design, then it is easily seen that for $i \neq j$ this rate of convergence is achieved by the estimator

$$\widehat{M}_{ij} = T^{-1} \sum_{k=1}^{T} (Y_{ik} - \bar{Y}_{\cdot k})(Y_{jk} - \bar{Y}_{\cdot k}), \qquad i \neq j,$$

and

$$\widehat{M}_{ii} = T^{-1} \sum_{k=1}^{T} (Y_{ik} - \bar{Y}_{\cdot k})^2 - \hat{\sigma}_i^2,$$

where $\hat{\sigma}_i^2$ denotes some nonparametric estimator of variance and $\bar{Y}_{\cdot k} = n^{-1} \times \sum_{j=1}^{n} Y_{jk}$.

In the case of a random design some adjustment is necessary: Define the ordered sample $t_{i(1)} \leq t_{i(2)} \leq \cdots \leq t_{i(T_i)}$ of design points, and for $j = 1, \ldots, T_i$, let $Y_{i(j)}$ denote the observation belonging to $t_{i(j)}$. With $t_{i(0)} = -t_{i(1)}$ and $t_{i(T_i+1)} = 2 - t_{i(T_i)}$, set

$$\chi_i(t) = \sum_{j=1}^{T_i} Y_{i(j)} I\left(t \in \left[\frac{t_{i(j-1)} + t_{i(j)}}{2}, \frac{t_{i(j)} + t_{i(j+1)}}{2}\right)\right), \qquad t \in [0, 1],$$

where $I(\cdot)$ denotes the indicator function, and for $i \neq j$, define the estimate of $M_{ij}$ by

$$\widehat{M}_{ij} = \int_0^1 \{\chi_i(t) - \bar{\chi}(t)\}\{\chi_j(t) - \bar{\chi}(t)\} \, dt,$$

where $\bar{\chi}(t) = n^{-1}\sum_{i=1}^{n}\chi_i(t)$. Finally, by redefining $t_{i(1)} = -t_{i(2)}$ and $t_{i(T_i+1)} = 2 - t_{i(T_i)}$, set $\chi_i^*(t) = \sum_{j=2}^{T_i} Y_{i(j-1)} I(t \in [\frac{t_{i(j-1)}+t_{i(j)}}{2}, \frac{t_{i(j)}+t_{i(j+1)}}{2}))$, $t \in [0,1]$. Then construct estimators of the diagonal terms $M_{ii}$ by

$$(7) \qquad \widehat{M}_{ii} = \int_0^1 \{\chi_i(t) - \bar{\chi}(t)\}\{\chi_i^*(t) - \bar{\chi}(t)\}\,dt.$$

The aim of using the estimator (7) for the diagonal terms is to avoid the additional bias implied by $\mathrm{E}_\varepsilon(Y_{ik}^2) = X_i(t_{ij})^2 + \sigma_i^2$. Here $\mathrm{E}_\varepsilon$ denotes conditional expectation given $t_{ij}$, $X_i$. Alternatively, we can construct a bias corrected estimator using some nonparametric estimation of variance $\sigma_i^2$, for example, the difference based model-free variance estimators studied in Hall, Kay and Titterington (1990) can be employed.

The eigenvalues $\hat{l}_1 \geq \hat{l}_2 \cdots$ and eigenvectors $\hat{p}_1, \hat{p}_2, \ldots$ of the resulting matrix $\widehat{M}$ then provide estimates $\hat{\lambda}_{r;T} = \hat{l}_r/n$ and $\hat{\beta}_{ri;T} = \sqrt{\hat{l}_r}\hat{p}_{ir}$ of $\hat{\lambda}_r$ and $\hat{\beta}_{ri}$. Estimates $\hat{\gamma}_{r;T}$ of the empirical functional principal component $\hat{\gamma}_r$ can be determined from (6) when replacing the unknown true functions $X_i$ by non-parametric estimates $\hat{X}_i$ (as, for example, local polynomial estimates) with smoothing parameter (bandwidth) $b$:

$$(8) \qquad \hat{\gamma}_{r;T} = \frac{1}{\sqrt{\hat{l}_r}}\sum_{i=1}^{n} \hat{p}_{ir}\hat{X}_i.$$

When considering (8), it is important to note that $\hat{\gamma}_{r;T}$ is defined as *a weighted average* of all estimated sample functions. Averaging reduces variance, and efficient estimation of $\hat{\gamma}_r$ therefore requires *undersmoothing* of individual function estimates $\hat{X}_i$. Theoretical results are given in Theorem 1 below. Indeed, if, for example, $n$ and $T = \min_i T_i$ are of the same order of magnitude, then under suitable additional regularity conditions it will be shown that for an optimal choice of a smoothing parameter $b \sim (nT)^{-1/5}$ and twice continuously differentiable $X_i$, we obtain the rate of convergence $\|\hat{\gamma}_r - \hat{\gamma}_{r;T}\| = \mathcal{O}_p\{(nT)^{-2/5}\}$. Note, however, that the bias corrected estimator (7) may yield negative eigenvalues. In practice, these values will be small and will have to be interpreted as zero. Furthermore, the eigenfunctions determined by (8) may not be exactly orthogonal. Again, when using reasonable bandwidths, this effect will be small, but of course (8) may by followed by a suitable orthogonalization procedure.

It is of interest to compare our procedure to more standard methods for estimating $\hat{\lambda}_r$ and $\hat{\gamma}_r$ as mentioned above. When evaluating eigenvalues and eigenfunctions of the empirical covariance operator of nonparametrically estimated curves $\hat{X}_i$, then for fixed $r \leq r_0$ the above rate of convergence for the estimated eigenfunctions may well be achieved for a suitable choice of smoothing parameters (e.g., number of basis functions). But as will be seen from Theorem 1, our approach also implies that $|\hat{\lambda}_r - \frac{\hat{l}_r}{n}| = \mathcal{O}_p(T^{-1} + n^{-1})$. When using standard methods it does not seem to be possible to obtain a corresponding rate of convergence, since any smoothing bias $|\mathrm{E}[\hat{X}_i(t)] - X_i(t)|$ will invariably affect the quality of the corresponding estimate of $\hat{\lambda}_r$.

We want to emphasize that any finite sample interpretation will require that $T$ is sufficiently large such that our nonparametric reconstructions of individual curves can be assumed to possess a fairly small bias. The above arguments do not apply to extremely sparse designs with very few observations per curve [see Hall, Müller and Wang (2006) for an FPCA methodology focusing on sparse data].

Note that, in addition to (8), our final estimate of the empirical mean function $\hat{\mu} = \bar{X}$ will be given by $\hat{\mu}_T = n^{-1} \sum_i \hat{X}_i$. A straightforward approach to determine a suitable bandwidth $b$ consists in a "leave-one-individual-out" cross-validation. For the maximal number $r_0$ of components to be considered, let $\hat{\mu}_{T,-i}$ and $\hat{\gamma}_{r;T,-i}, r = 1, \ldots, r_0$, denote the estimates of $\hat{\mu}$ and $\hat{\gamma}_r$ obtained from the data $(Y_{lj}, t_{lj}), l = 1, \ldots, i-1, i+1, \ldots, n, j = 1, \ldots, T_k$. By (8), these estimates depend on $b$, and one may approximate an optimal smoothing parameter by minimizing

$$\sum_i \sum_j \left\{ Y_{ij} - \hat{\mu}_{T,-i}(t_{ij}) - \sum_{r=1}^{r_0} \hat{\vartheta}_{ri} \hat{\gamma}_{r;T,-i}(t_{ij}) \right\}^2$$

over $b$, where $\hat{\vartheta}_{ri}$ denote ordinary least squares estimates of $\hat{\beta}_{ri}$. A more sophisticated version of this method may even allow to select different bandwidths $b_r$ when estimating different functional principal components by (8). Although, under regularity conditions, the same qualitative rates of convergence hold for any arbitrary *fixed* $r \le r_0$, the quality of estimates decreases when $r$ becomes large. Due to $\langle \gamma_s, \gamma_r \rangle = 0$ for $s < r$, the number of zero crossings, peaks and valleys of $\gamma_r$ has to increase with $r$. Hence, in tendency $\gamma_r$ will be less and less smooth as $r$ increases. At the same time, $\lambda_r \to 0$, which means that for large $r$ the $r$th eigenfunctions will only possess a very small influence on the structure of $X_i$. This in turn means that the relative importance of the error terms $\varepsilon_{ik}$ in (4) on the structure of $\hat{\gamma}_{r;T}$ will increase with $r$.

### 2.2. *One sample inference.*

Clearly, in the framework described by (1)–(4) we are faced with two sources of variability of estimated functional principal components. Due to sampling variation, $\hat{\gamma}_r$ will differ from the true component $\gamma_r$, and due to (4), there will exist an additional estimation error when approximating $\hat{\gamma}_r$ by $\hat{\gamma}_{r;T}$.

The following theorems quantify the order of magnitude of these different types of error. Our theoretical results are based on the following assumptions on the structure of the random functions $X_i$.

ASSUMPTION 1.    $X_1, \ldots, X_n \in L^2[0,1]$ is an i.i.d. sample of random functions with mean $\mu$ and continuous covariance function $\sigma(t,s)$, and (1) holds for a system of eigenfunctions satisfying $\sup_{s \in \mathbb{N}} \sup_{t \in [0,1]} |\gamma_s(t)| < \infty$. Furthermore, $\sum_{r=1}^{\infty} \sum_{s=1}^{\infty} \mathrm{E}[\beta_{ri}^2 \beta_{si}^2] < \infty$ and $\sum_{q=1}^{\infty} \sum_{s=1}^{\infty} \mathrm{E}[\beta_{ri}^2 \beta_{qi} \beta_{si}] < \infty$ for all $r \in \mathbb{N}$.

Recall that $\mathrm{E}[\beta_{ri}] = 0$ and $\mathrm{E}[\beta_{ri} \beta_{si}] = 0$ for $r \neq s$. Note that the assumption on the factor loadings is necessarily fulfilled if $X_i$ are Gaussian random functions. Then $\beta_{ri}$ and $\beta_{si}$ are independent for $r \neq s$, all moments of $\beta_{ri}$ are finite, and hence $\mathrm{E}[\beta_{ri}^2 \beta_{qi} \beta_{si}] = 0$ for $q \neq s$, as well as $\mathrm{E}[\beta_{ri}^2 \beta_{si}^2] = \lambda_r \lambda_s$ for $r \neq s$; see Gihman and Skorohod (1973).

We need some further assumptions concerning smoothness of $X_i$ and the structure of the discrete model (4).

ASSUMPTION 2.    (a) $X_i$ is a.s. twice continuously differentiable. There exists a constant $D_1 < \infty$ such that the derivatives are bounded by $\sup_t \mathrm{E}[X_i'(t)^4] \leq D_1$, as well as $\sup_t \mathrm{E}[X_i''(t)^4] \leq D_1$.

(b) The design points $t_{ik}$, $i = 1, \ldots, n$, $k = 1, \ldots, T_i$, are i.i.d. random variables which are independent of $X_i$ and $\varepsilon_{ik}$. The corresponding design density $f$ is continuous on $[0,1]$ and satisfies $\inf_{t \in [0,1]} f(t) > 0$.

(c) For any $i$, the error terms $\varepsilon_{ik}$ are i.i.d. zero mean random variables with $\mathrm{Var}(\varepsilon_{ik}) = \sigma_i^2$. Furthermore, $\varepsilon_{ik}$ is independent of $X_i$, and there exists a constant $D_2$ such that $\mathrm{E}(\varepsilon_{ik}^8) < D_2$ for all $i, k$.

(d) The estimates $\hat{X}_i$ used in (8) are determined by either a local linear or a Nadaraya–Watson kernel estimator with smoothing parameter $b$ and kernel function $K$. $K$ is a continuous probability density which is symmetric at 0.

The following theorems provide asymptotic results as $n, T \to \infty$, where $T = \min_{i=1}^{n} \{T_i\}$.

THEOREM 1.    *In addition to Assumptions 1 and 2, assume that* $\inf_{s \neq r} |\lambda_r - \lambda_s| > 0$ *holds for some* $r = 1, 2, \ldots$ *. Then we have the following:*

(i) $n^{-1} \sum_{i=1}^{n} (\hat{\beta}_{ri} - \hat{\beta}_{ri;T})^2 = \mathcal{O}_p(T^{-1})$ *and*

$$(9) \qquad \left| \hat{\lambda}_r - \frac{\hat{l}_r}{n} \right| = \mathcal{O}_p(T^{-1} + n^{-1}).$$

(ii) *If additionally* $b \to 0$ *and* $(Tb)^{-1} \to 0$ *as* $n, T \to \infty$, *then for all* $t \in [0,1]$,

$$(10) \qquad |\hat{\gamma}_r(t) - \hat{\gamma}_{r;T}(t)| = \mathcal{O}_p\{b^2 + (nTb)^{-1/2} + (Tb^{1/2})^{-1} + n^{-1}\}.$$

*A proof is given in the Appendix.*

THEOREM 2. *Under Assumption 1 we obtain the following:*

(i) *For all $t \in [0, 1]$,*

$$\sqrt{n}\{\bar{X}(t) - \mu(t)\} = \sum_r \left\{ \frac{1}{\sqrt{n}} \sum_{i=1}^n \beta_{ri} \right\} \gamma_r(t) \xrightarrow{\mathcal{L}} N\left(0, \sum_r \lambda_r \gamma_r(t)^2\right).$$

*If, furthermore, $\lambda_{r-1} > \lambda_r > \lambda_{r+1}$ holds for some fixed $r \in \{1, 2, \ldots\}$, then*

(ii)

$$(11) \qquad \sqrt{n}(\hat{\lambda}_r - \lambda_r) = \frac{1}{\sqrt{n}} \sum_{i=1}^n (\beta_{ri}^2 - \lambda_r) + \mathcal{O}_p(n^{-1/2}) \xrightarrow{\mathcal{L}} N(0, \Lambda_r),$$

*where $\Lambda_r = \mathrm{E}[(\beta_{ri}^2 - \lambda_r)^2]$,*

(iii) *and for all $t \in [0, 1]$*

$$(12) \qquad \hat{\gamma}_r(t) - \gamma_r(t) = \sum_{s \neq r} \left\{ \frac{1}{n(\lambda_r - \lambda_s)} \sum_{i=1}^n \beta_{si} \beta_{ri} \right\} \gamma_s(t) + R_r(t),$$

$$where \ \|R_r\| = \mathcal{O}_p(n^{-1}).$$

*Moreover,*

$$\sqrt{n} \sum_{s \neq r} \left\{ \frac{1}{n(\lambda_r - \lambda_s)} \sum_{i=1}^n \beta_{si} \beta_{ri} \right\} \gamma_s(t)$$

$$\xrightarrow{\mathcal{L}} N\left(0, \sum_{q \neq r} \sum_{s \neq r} \frac{\mathrm{E}[\beta_{ri}^2 \beta_{qi} \beta_{si}]}{(\lambda_q - \lambda_r)(\lambda_s - \lambda_r)} \gamma_q(t) \gamma_s(t)\right).$$

A proof can be found in the Appendix. The theorem provides a generalization of the results of Dauxois, Pousse and Romain (1982) who derive explicit asymptotic distributions by assuming Gaussian random functions $X_i$. Note that in this case $\Lambda_r = 2\lambda_r^2$ and $\sum_{q \neq r} \sum_{s \neq r} \frac{\mathrm{E}[\beta_{ri}^2 \beta_{qi} \beta_{si}]}{(\lambda_q - \lambda_r)(\lambda_s - \lambda_r)} \gamma_q(t) \gamma_s(t) = \sum_{s \neq r} \frac{\lambda_r \lambda_s}{(\lambda_s - \lambda_r)^2} \gamma_s(t)^2$.

When evaluating the bandwidth-dependent terms in (10), best rates of convergence $|\hat{\gamma}_r(t) - \hat{\gamma}_{r;T}(t)| = \mathcal{O}_p\{(nT)^{-2/5} + T^{-4/5} + n^{-1}\}$ are achieved when choosing an undersmoothing bandwidth $b \sim \max\{(nT)^{-1/5}, T^{-2/5}\}$. Theoretical work in functional data analysis is usually based on the implicit assumption that the additional error due to (4) is negligible, and that one can proceed "as if" the functions $X_i$ were directly observed. In view of Theorems 1 and 2, this approach is justified in the following situations:

(1) $T$ is much larger than $n$, that is, $n/T^{4/5} \to 0$, and the smoothing parameter $b$ in (8) is of order $T^{-1/5}$ (optimal smoothing of individual functions).

(2) An undersmoothing bandwidth $b \sim \max\{(nT)^{-1/5}, T^{-2/5}\}$ is used and $n/T^{8/5} \to 0$. This means that $T$ may be smaller than $n$, but $T$ must be at least of order of magnitude larger than $n^{5/8}$.

In both cases (1) and (2) the above theorems imply that $|\hat{\lambda}_r - \frac{\tilde{l}_r}{n}| = o_p(|\hat{\lambda}_r - \lambda_r|)$, as well as $\|\hat{\gamma}_r - \hat{\gamma}_{r;T}\| = o_p(\|\hat{\gamma}_r - \gamma_r\|)$. Inference about functional principal components will then be first-order equivalent to an inference based on known functions $X_i$.

In such situations Theorem 2 suggests bootstrap procedures as tools for one sample inference. For example, the distribution of $\|\hat{\gamma}_r - \gamma_r\|$ may by approximated by the bootstrap distribution of $\|\hat{\gamma}_r^* - \hat{\gamma}_r\|$, where $\hat{\gamma}_r^*$ are estimates to be obtained from i.i.d. bootstrap resamples $X_1^*, X_2^*, \ldots, X_n^*$ of $\{X_1, X_2, \ldots, X_n\}$. This means that $X_1^* = X_{i_1}, \ldots, X_n^* = X_{i_n}$ for some indices $i_1, \ldots, i_n$ drawn independently and with replacement from $\{1, \ldots, n\}$ and, in practice, $\hat{\gamma}_r^*$ may thus be approximated from corresponding discrete data $(Y_{i_1 j}, t_{i_1 j})_{j=1,\ldots,T_{i_1}}, \ldots, (Y_{i_n j}, t_{i_n j})_{j=1,\ldots,T_{i_n}}$. The additional error is negligible if either (1) or (2) is satisfied.

One may wonder about the validity of such a bootstrap. Functions are complex objects and there is no established result in bootstrap theory which readily generalizes to samples of random functions. But by (1), i.i.d. bootstrap resamples $\{X_i^*\}_{i=1,\ldots,n}$ may be equivalently represented by corresponding, i.i.d. resamples $\{\beta_{1i}^*, \beta_{2i}^*, \ldots\}_{i=1,\ldots,n}$ of factor loadings. Standard multivariate bootstrap theorems imply that for any $q \in \mathbb{N}$ the distribution of moments of the random vectors $(\beta_{1i}, \ldots, \beta_{qi})$ may be consistently approximated by the bootstrap distribution of corresponding moments of $(\beta_{1i}^*, \ldots, \beta_{qi}^*)$. Together with some straightforward limit arguments as $q \to \infty$, the structure of the first-order terms in the asymptotic expansions (11) and (12) then allows to establish consistency of the functional bootstrap. These arguments will be made precise in the proof of Theorem 3 below, which concerns related bootstrap statistics in two sample problems.

REMARK. Theorem 2(iii) implies that the variance of $\hat{\gamma}_r$ is large if one of the differences $\lambda_{r-1} - \lambda_r$ or $\lambda_r - \lambda_{r+1}$ is small. In the limit case of eigenvalues of multiplicity $m > 1$ our theory does not apply. Note that then only the $m$-dimensional eigenspace is identified, but not a particular basis (eigenfunctions). In multivariate PCA Tyler (1981) provides some inference results on corresponding projection matrices assuming that $\lambda_r > \lambda_{r+1} \geq \cdots \geq \lambda_{r+m} > \lambda_{r+m+1}$ for known values of $r$ and $m$.

Although the existence of eigenvalues $\lambda_r$, $r \leq r_0$, with multiplicity $m > 1$ may be considered as a degenerate case, it is immediately seen that $\lambda_r \to 0$ and, hence, $\lambda_r - \lambda_{r+1} \to 0$ as $r$ increases. Even in the case of fully observed functions $X_i$, estimates of eigenfunctions corresponding to very small eigenvalues will thus be poor. The problem of determining a sensible upper limit of the number $r_0$ of principal components to be analyzed is addressed in Hall and Hosseini-Nasab (2006).

**3. Two sample inference.** The comparison of functional components across groups leads naturally to two sample problems. Thus, let

$$X_1^{(1)}, X_2^{(1)}, \ldots, X_{n_1}^{(1)} \quad \text{and} \quad X_1^{(2)}, X_2^{(2)}, \ldots, X_{n_2}^{(2)}$$

denote two independent samples of smooth functions. The problem of interest is to test in how far the distributions of these random functions coincide. The structure of the different distributions in function space can be accessed by means of the respective Karhunen–Loève decompositions. The problem to be considered then translates into testing equality of the different components of these decompositions given by

$$(13) \qquad X_i^{(p)} = \mu^{(p)} + \sum_{r=1}^{\infty} \beta_{ri}^{(p)} \gamma_r^{(p)}, \qquad p = 1, 2,$$

where again $\gamma_r^{(p)}$ are the eigenfunctions of the respective covariance operator $\Gamma^{(p)}$ corresponding to the eigenvalues $\lambda_1^{(p)} = \mathrm{E}\{(\beta_{1i}^{(p)})^2\} \geq \lambda_2^{(p)} = \mathrm{E}\{(\beta_{2i}^{(p)})^2\} \geq \cdots$. We will again suppose that $\lambda_{r-1}^{(p)} > \lambda_r^{(p)} > \lambda_{r+1}^{(p)}$, $p = 1, 2$, for all $r \leq r_0$ components to be considered. Without restriction, we will additionally assume that signs are such that $\langle \gamma_r^{(1)}, \gamma_r^{(2)} \rangle \geq 0$, as well as $\langle \hat{\gamma}_r^{(1)}, \hat{\gamma}_r^{(2)} \rangle \geq 0$.

It is of great interest to detect possible variations in the functional components characterizing the two samples in (13). Significant difference may give rise to substantial interpretation. Important hypotheses to be considered thus are as follows:

$$H_{0_1}: \mu^{(1)} = \mu^{(2)} \quad \text{and} \quad H_{0_{2,r}}: \gamma_r^{(1)} = \gamma_r^{(2)}, \qquad r \leq r_0.$$

Hypothesis $H_{0_{2,r}}$ is of particular importance. Then $\gamma_r^{(1)} = \gamma_r^{(2)}$ and only the factor loadings $\beta_{ri}$ may vary across samples. If, for example, $H_{0_{2,r}}$ is accepted, one may additionally want to test hypotheses about the distributions of $\beta_{ri}^{(p)}$, $p = 1, 2$. Recall that necessarily $\mathrm{E}\{\beta_{ri}^{(p)}\} = 0$, $\mathrm{E}\{\beta_{ri}^{(p)}\}^2 = \lambda_r^{(p)}$, and $\beta_{si}^{(p)}$ is uncorrelated with $\beta_{ri}^{(p)}$ if $r \neq s$. If the $X_i^{(p)}$ are Gaussian random variables, the $\beta_{ri}^{(p)}$ are independent $N(0, \lambda_r)$ random variables. A natural hypothesis to be tested then refers to the equality of variances:

$$H_{0_{3,r}}: \lambda_r^{(1)} = \lambda_r^{(2)}, \qquad r = 1, 2, \ldots.$$

Let $\hat{\mu}^{(p)}(t) = \frac{1}{n_p} \sum_i X_i^{(p)}(t)$, and let $\hat{\lambda}_1^{(p)} \geq \hat{\lambda}_2^{(p)} \geq \cdots$ and $\hat{\gamma}_1^{(p)}, \hat{\gamma}_2^{(p)}, \ldots$ denote eigenvalues and corresponding eigenfunctions of the empirical covariance operator $\hat{\Gamma}_{n_p}^{(p)}$ of $X_1^{(p)}, X_2^{(p)}(t), \ldots, X_{n_p}^{(p)}$. The following test statistics are defined in terms of $\hat{\mu}^{(p)}$, $\hat{\lambda}_r^{(p)}$ and $\hat{\gamma}_r^{(p)}$. As discussed in the proceeding section, all curves in both samples are usually not directly observed, but have to be reconstructed from noisy observations according to (4). In this situation, the "true" empirical eigenvalues and eigenfunctions have to be replaced by their discrete sample estimates. Bootstrap estimates are obtained by resampling the observations corresponding to the unknown curves $X_i^{(p)}$. As discussed in Section 2.2, the validity of our test procedures is then based on the assumption that $T$ is sufficiently large such that the additional estimation error is asymptotically negligible.

Our tests of the hypotheses $H_{0_1}, H_{0_{2,r}}$ and $H_{0_{3,r}}$ rely on the statistics

$$D_1 \stackrel{\text{def}}{=} \|\hat{\mu}^{(1)} - \hat{\mu}^{(2)}\|^2,$$

$$D_{2,r} \stackrel{\text{def}}{=} \|\hat{\gamma}_r^{(1)} - \hat{\gamma}_r^{(2)}\|^2,$$

$$D_{3,r} \stackrel{\text{def}}{=} |\hat{\lambda}_r^{(1)} - \hat{\lambda}_r^{(2)}|^2.$$

The respective null-hypothesis has to be rejected if $D_1 \geq \Delta_{1;1-\alpha}$, $D_{2,r} \geq \Delta_{2,r;1-\alpha}$ or $D_{3,r} \geq \Delta_{3,r;1-\alpha}$, where $\Delta_{1;1-\alpha}$, $\Delta_{2,r;1-\alpha}$ and $\Delta_{3,r;1-\alpha}$ denote the critical values of the distributions of

$$\Delta_1 \stackrel{\text{def}}{=} \|\hat{\mu}^{(1)} - \mu^{(1)} - (\hat{\mu}^{(2)} - \mu^{(2)})\|^2,$$

$$\Delta_{2,r} \stackrel{\text{def}}{=} \|\hat{\gamma}_r^{(1)} - \gamma_r^{(1)} - (\hat{\gamma}_r^{(2)} - \gamma_r^{(2)})\|^2,$$

$$\Delta_{3,r} \stackrel{\text{def}}{=} |\hat{\lambda}_r^{(1)} - \lambda_r^{(1)} - (\hat{\lambda}_r^{(2)} - \lambda_r^{(2)})|^2.$$

Of course, the distributions of the different $\Delta$'s cannot be accessed directly, since they depend on the unknown true population mean, eigenvalues and eigenfunctions. However, it will be shown below that these distributions and, hence, their critical values are approximated by the bootstrap distribution of

$$\Delta_1^* \stackrel{\text{def}}{=} \|\hat{\mu}^{(1)*} - \hat{\mu}^{(1)} - (\hat{\mu}^{(2)*} - \hat{\mu}^{(2)})\|^2,$$

$$\Delta_{2,r}^* \stackrel{\text{def}}{=} \|\hat{\gamma}_r^{(1)*} - \hat{\gamma}_r^{(1)} - (\hat{\gamma}_r^{(2)*} - \hat{\gamma}_r^{(2)})\|^2,$$

$$\Delta_{3,r}^* \stackrel{\text{def}}{=} |\hat{\lambda}_r^{(1)*} - \hat{\lambda}_r^{(1)} - (\hat{\lambda}_r^{(2)*} - \hat{\lambda}_r^{(2)})|^2,$$

where $\hat{\mu}^{(1)*}$, $\hat{\gamma}_r^{(1)*}$, $\hat{\lambda}_r^{(1)*}$, as well as $\hat{\mu}^{(2)*}$, $\hat{\gamma}_r^{(2)*}$, $\hat{\lambda}_r^{(2)*}$, are estimates to be obtained from independent bootstrap samples $X_1^{1*}(t), X_2^{1*}(t), \ldots, X_{n_1}^{1*}(t)$, as well as $X_1^{2*}(t), X_2^{2*}(t), \ldots, X_{n_2}^{2*}(t)$.

This test procedure is motivated by the following insights:

(1) Under each of our null-hypotheses the respective test statistics $D$ is equal to the corresponding $\Delta$. The test will thus asymptotically possess the correct level: $P(D > \Delta_{1-\alpha}) \approx \alpha$.

(2) If the null hypothesis is false, then $D \neq \Delta$. Compared to the distribution of $\Delta$, the distribution of $D$ is shifted by the difference in the true means, eigenfunctions or eigenvalues. In tendency $D$ will be larger than $\Delta_{1-\alpha}$.

Let $1 < L \leq r_0$. Even if for $r \leq L$ the equality of eigenfunctions is rejected, we may be interested in the question of whether at least the $L$-dimensional eigenspaces generated by the first $L$ eigenfunctions are identical. Therefore, let $\mathcal{E}_L^{(1)}$, as well as $\mathcal{E}_L^{(2)}$, denote the $L$-dimensional linear function spaces generated by the eigenfunctions $\gamma_1^{(1)}, \ldots, \gamma_L^{(1)}$ and $\gamma_1^{(2)}, \ldots, \gamma_L^{(2)}$, respectively. We then aim to test the null hypothesis:

$$H_{0_{4,L}} : \mathcal{E}_L^{(1)} = \mathcal{E}_L^{(2)}.$$

Of course, $H_{0_{4,L}}$ corresponds to the hypothesis that the operators projecting into $\mathcal{E}_L^{(1)}$ and $\mathcal{E}_L^{(2)}$ are identical. This in turn translates into the condition that

$$\sum_{r=1}^{L} \gamma_r^{(1)}(t)\gamma_r^{(1)}(s) = \sum_{r=1}^{L} \gamma_r^{(2)}(t)\gamma_r^{(2)}(s) \qquad \text{for all } t, s \in [0, 1].$$

Similar to above, a suitable test statistic is given by

$$D_{4,L} \stackrel{\text{def}}{=} \iint \left\{ \sum_{r=1}^{L} \hat{\gamma}_r^{(1)}(t)\hat{\gamma}_r^{(1)}(s) - \sum_{r=1}^{L} \hat{\gamma}_r^{(2)}(t)\hat{\gamma}_r^{(2)}(s) \right\}^2 dt\, ds$$

and the null hypothesis is rejected if $D_{4,L} \geq \Delta_{4,L;1-\alpha}$, where $\Delta_{4,L;1-\alpha}$ denotes the critical value of the distribution of

$$\Delta_{4,L} \stackrel{\text{def}}{=} \iint \left[ \sum_{r=1}^{L} \{\hat{\gamma}_r^{(1)}(t)\hat{\gamma}_r^{(1)}(s) - \gamma_r^{(1)}(t)\gamma_r^{(1)}(s)\} \right.$$

$$\left. - \sum_{r=1}^{L} \{\hat{\gamma}_r^{(2)}(t)\hat{\gamma}_r^{(2)}(s) - \gamma_r^{(2)}(t)\gamma_r^{(2)}(s)\} \right]^2 dt\, ds.$$

The distribution of $\Delta_{4,L}$ and, hence, its critical values are approximated by the bootstrap distribution of

$$\Delta_{4,L}^* \stackrel{\text{def}}{=} \iint \left[ \sum_{r=1}^{L} \{\hat{\gamma}_r^{(1)*}(t)\hat{\gamma}_r^{(1)*}(s) - \hat{\gamma}_r^{(1)}(t)\hat{\gamma}_r^{(1)}(s)\} \right.$$

$$\left. - \sum_{r=1}^{L} \{\hat{\gamma}_r^{(2)*}(t)\hat{\gamma}_r^{(2)*}(s) - \hat{\gamma}_r^{(2)}(t)\hat{\gamma}_r^{(2)}(s)\} \right]^2 dt\, ds.$$

It will be shown in Theorem 3 below that under the null hypothesis, as well as under the alternative, the distributions of $n\Delta_1, n\Delta_{2,r}, n\Delta_{3,r}, n\Delta_{4,L}$ converge to continuous limit distributions which can be consistently approximated by the bootstrap distributions of $n\Delta_1^*, n\Delta_{2,r}^*, n\Delta_{3,r}^*, n\Delta_{4,L}^*$.

3.1. *Theoretical results.* Let $n = (n_1 + n_2)/2$. We will assume that asymptotically $n_1 = n \cdot q_1$ and $n_2 = n \cdot q_2$ for some fixed proportions $q_1$ and $q_2$. We will then study the asymptotic behavior of our statistics as $n \to \infty$.

We will use $\mathcal{X}_1 = \{X_1^{(1)}, \ldots, X_{n_1}^{(1)}\}$ and $\mathcal{X}_2 = \{X_1^{(2)}, \ldots, X_{n_2}^{(2)}\}$ to denote the observed samples of random functions.

THEOREM 3. *Assume that* $\{X_1^{(1)}, \ldots, X_{n_1}^{(1)}\}$ *and* $\{X_1^{(2)}, \ldots, X_{n_2}^{(2)}\}$ *are two independent samples of random functions, each of which satisfies Assumption 1. As $n \to \infty$ we then obtain the following:*

(i) *There exists a nondegenerated, continuous probability distribution $F_1$ such that $n\Delta_1 \xrightarrow{\mathcal{L}} F_1$, and for any $\delta > 0$,*

$$|P(n\Delta_1 \geq \delta) - P(n\Delta_1^* \geq \delta | \mathcal{X}_1, \mathcal{X}_2)| = \mathcal{o}_p(1).$$

(ii) *If, furthermore, $\lambda_{r-1}^{(1)} > \lambda_r^{(1)} > \lambda_{r+1}^{(1)}$ and $\lambda_{r-1}^{(2)} > \lambda_r^{(2)} > \lambda_{r+1}^{(2)}$ hold for some fixed $r = 1, 2, \ldots$, there exist a nondegenerated, continuous probability distributions $F_{k,r}$ such that $n\Delta_{k,r} \xrightarrow{\mathcal{L}} F_{k,r}$, $k = 2, 3$, and for any $\delta > 0$,*

$$|P(n\Delta_{k,r} \geq \delta) - P(n\Delta_{k,r}^* \geq \delta | \mathcal{X}_1, \mathcal{X}_2)| = \mathcal{o}_p(1), \qquad k = 2, 3.$$

(iii) *If $\lambda_r^{(1)} > \lambda_{r+1}^{(1)} > 0$ and $\lambda_r^{(2)} > \lambda_{r+1}^{(2)} > 0$ hold for all $r = 1, \ldots, L$, there exists a nondegenerated, continuous probability distribution $F_{4,L}$ such that $n\Delta_{4,L} \xrightarrow{\mathcal{L}} F_{4,L}$, and for any $\delta > 0$,*

$$|P(n\Delta_{4,L} \geq \delta) - P(n\Delta_{4,L}^* \geq \delta | \mathcal{X}_1, \mathcal{X}_2)| = \mathcal{o}_p(1).$$

The structures of the distributions $F_1$, $F_{2,r}$, $F_{3,r}$, $F_{4,L}$ are derived in the proof of the theorem which can be found in the Appendix. They are obtained as limits of distributions of quadratic forms.

3.2. *Simulation study.* In this paragraph we illustrate the finite behavior of the proposed test. The basic simulation-setup (setup "a") is established as follows: the first sample is generated by the random combination of orthonormalized sine and cosine functions (Fourier functions) and the second sample is generated by the random combination of the same but shifted factor functions:

$$X_i^{(1)}(t_{ik}) = \beta_{1i}^{(1)} \sqrt{2} \sin(2\pi t_{ik}) + \beta_{2i}^{(1)} \sqrt{2} \cos(2\pi t_{ik}),$$

$$X_i^{(2)}(t_{ik}) = \beta_{1i}^{(2)} \sqrt{2} \sin\{2\pi(t_{ik} + \delta)\} + \beta_{2i}^{(2)} \sqrt{2} \cos\{2\pi(t_{ik} + \delta)\}.$$

The factor loadings are i.i.d. random variables with $\beta_{1i}^{(p)} \sim N(0, \lambda_1^{(p)})$ and $\beta_{2i}^{(p)} \sim N(0, \lambda_2^{(p)})$. The functions are generated on the equidistant grid $t_{ik} = t_k = k/T$, $k = 1, \ldots T = 100$, $i = 1, \ldots, n = 70$. The simulation setup is based

| Setup/shift | 0 | 0.05 | 0.1 | 0.15 | 0.2 | 0.25 |
|---|---|---|---|---|---|---|
| (a) 10, 5, 8, 4 | 0.13 | 0.41 | 0.85 | 0.96 | 1 | 1 |
| (a) 4, 2, 2, 1 | 0.12 | 0.48 | 0.87 | 0.96 | 1 | 1 |
| (a) 2, 1, 1.5, 2 | 0.14 | 0.372 | 0.704 | 0.872 | 0.92 | 0.9 |
| (b) 10, 5, 8, 4 $D_1$ | 0.10 | 0.44 | 0.86 | 0.95 | 1 | 1 |
| (b) 10, 5, 8, 4 $D_2$ | 1 | 1 | 1 | 1 | 1 | 1 |

on the fact that the error of the estimation of the eigenfunctions simulated by sine and cosine functions is, in particular, manifested by some shift of the estimated eigenfunctions. The focus of this simulation study is the test of common eigenfunctions.

For the presentation of results in Table 1, we use the following notation: "(a) $\lambda_1^{(1)}, \lambda_2^{(1)}, \lambda_2^{(2)}, \lambda_2^{(2)}$." The shift parameter $\delta$ is changing from 0 to 0.25 with the step 0.05. It should be mentioned that the shift $\delta = 0$ yields the simulation of level and setup with shift $\delta = 0.25$ yields the simulation of the alternative, where the two factor functions are exchanged.

In the second setup (setup "b") the first factor functions are the same and the second factor functions differ:

$$X_i^{(1)}(t_{ik}) = \beta_{1i}^{(1)}\sqrt{2}\sin(2\pi t_{ik}) + \beta_{2i}^{(1)}\sqrt{2}\cos(2\pi t_{ik}),$$

$$X_i^{(2)}(t_{ik}) = \beta_{1i}^{(2)}\sqrt{2}\sin\{2\pi(t_{ik}+\delta)\} + \beta_{2i}^{(2)}\sqrt{2}\sin\{4\pi(t_{ik}+\delta)\}.$$

In Table 1 we use the notation "(b) $\lambda_1^{(1)}, \lambda_2^{(1)}, \lambda_2^{(2)}, \lambda_2^{(2)}, D_r$." $D_r$ means the test for the equality of the $r$th eigenfunction. In the bootstrap tests we used 500 bootstrap replications. The critical level in this simulation is $\alpha = 0.1$. The number of simulations is 250.

We can interpret Table 1 in the following way: In power simulations ($\delta \neq 0$) test behaves as expected: less powerful if the functions are "hardly distinguishable" (small shift, small difference in eigenvalues). The level approximation seems to be less precise if the difference in the eingenvalues $(\lambda_1^{(p)} - \lambda_2^{(p)})$ becomes smaller. This can be explained by relative small sample-size $n$, small number of bootstrap-replications and increasing estimation-error as argued in Theorem 2, assertion (iii).

In comparison to our general setup (4), we used an equidistant and common design for all functions. This simplification is necessary, it simplifies and speeds-up the simulations, in particular, using general random and observation-specific design makes the simulation computationally untractable.

Second, we omitted the additional observation error, this corresponds to the standard assumptions in the functional principal components theory. As argued in Section 2.2, the inference based on the directly observed functions and estimated functions $X_i$ is first-order equivalent under mild conditions implied by Theorems 1 and 2. In order to illustrate this theoretical result in the simulation, we used the following setup:

$$X_i^{(1)}(t_{ik}) = \beta_{1i}^{(1)} \sqrt{2} \sin(2\pi t_{ik}) + \beta_{2i}^{(1)} \sqrt{2} \cos(2\pi t_{ik}) + \varepsilon_{ik}^{(1)},$$

$$X_i^{(2)}(t_{ik}) = \beta_{1i}^{(2)} \sqrt{2} \sin\{2\pi(t_{ik} + \delta)\} + \beta_{2i}^{(2)} \sqrt{2} \cos\{2\pi(t_{ik} + \delta)\} + \varepsilon_{ik}^{(2)},$$

where $\varepsilon_{ik}^{(p)} \sim N(0, 0.25)$, $p = 1, 2$, all other parameters remain the same as in the simulation setup "a." Using this setup, we recalculate the simulation presented in the second "row" of Table 1, for estimation of the functions $X_i^{(p)}, p = 1, 2$, we used the Nadaraya–Watson estimation with Epanechnikov kernel and bandwidth $b = 0.05$. We run the simulations with various bandwidths, the choice of the bandwidth does not have a strong influence on results except by oversmoothing (large bandwidths). The results are printed in Table 2. As we can see, the difference of the simulation results using estimated functions is not significant in comparison to the results printed in the second line of Table 1—directly observed functional values.

Table 2
*The results of the simulation for $\alpha = 0.1$, $n = 70$, $T = 100$ with additional error in observation*

| Setup/shift | 0 | 0.05 | 0.1 | 0.15 | 0.2 | 0.25 |
|---|---|---|---|---|---|---|
| (a) 10, 5, 8, 4 | 0.09 | 0.35 | 0.64 | 0.92 | 0.94 | 0.97 |

The last limitation of this simulation study is the choice of a particular alternative. A more general setup of this simulation study might be based on the following model: $X_i^{(1)}(t) = \beta_{1i}^{(1)} \gamma_1^{(1)}(t) + \beta_{2i}^{(1)} \gamma_2^{(1)}(t)$, $X_i^{(2)}(t) = \beta_{1i}^{(2)} \gamma_1^{(2)}(t) + \beta_{2i}^{(2)} \gamma_2^{(2)}(t)$, where $\gamma_1^{(1)}, \gamma_1^{(2)}, \gamma_2^{(2)}$ and $g$ are mutually orthogonal functions on $L^2[0, 1]$ and $\gamma_2^{(2)} = (1 + v^2)^{-1/2}\{\gamma_2^{(1)} + vg\}$. Basically we create the alternative by the contamination of one of the "eigenfunctions" (in our case the second one) in the direction $g$ and ensure $\|\gamma_2^{(2)}\| = 1$. The amount of the contamination is controlled by the parameter $v$. Note that the exact squared integral difference $\|\gamma_2^{(1)} - \gamma_2^{(2)}\|^2$ does not depend on function $g$. Thus, in the "functional sense" particular "direction of the alternative hypothesis" represented by the function $g$ has no impact on the power of the test. However, since we are using a nonparametric estimation technique, we might expect that rough (highly fluctuating) functions $g$ will yield higher error of estimation and, hence, decrease the precision (and power) of the test. Finally, a higher number of factor functions (L) in simulation may cause less precise approximation of critical values and more bootstrap replications and larger sample-size may be needed. This can also be expected from Theorem 2 in Section 2.2—the variance of the estimated eigenfunctions depends on all eigenfunctions corresponding to nonzero eigenvalues.

**4. Implied volatility analysis.** In this section we present an application of the method discussed in previous sections to the implied volatilities of European options on the German stock index (ODAX). Implied volatilities are derived from the Black–Scholes (BS) pricing formula for European options; see Black and Scholes (1973). European call and put options are derivatives written on an underlying asset with price process $S_i$, which yield the pay-off $\max(S_I - K, 0)$ and $\max(K - S_I, 0)$, respectively. Here $i$ denotes the current day, $I$ the expiration day and $K$ the strike price. Time to maturity is defined as $\tau = I - i$. The BS pricing formula for a Call option is

$$(14) \qquad C_i(S_i, K, \tau, r, \sigma) = S_i \Phi(d_1) - K e^{-r\tau} \Phi(d_2),$$

where $d_1 = \frac{\ln(S_i/K) + (r + \sigma^2/2)\tau}{\sigma\sqrt{\tau}}$, $d_2 = d_1 - \sigma\sqrt{\tau}$, $r$ is the risk-free interest rate, $\sigma$ is the (unknown and constant) volatility parameter, and $\Phi$ denotes the c.d.f. of a standard normal distributed random variable. In (14) we assume the zero-dividend case. The Put option price $P_i$ can be obtained from the put–call parity $P_i = C_i - S_i + e^{-r\tau} K$.

The implied volatility $\tilde{\sigma}$ is defined as the volatility $\sigma$, for which the BS price $C_i$ in (14) equals the price $\tilde{C}_i$ observed on the market. For a single asset, we obtain at each time point (day $i$) and for each maturity $\tau$ a IV function $\tilde{\sigma}_i^\tau(K)$. Practitioners often rescale the strike dimension by plotting this surface in terms of (futures) moneyness $\kappa = K/F_i(\tau)$, where $F_i(\tau) = S_i e^{r\tau}$.

Clearly, for given parameters $S_i, r, K, \tau$ the mapping from prices to IVs is a one-to-one mapping. The IV is often used for quoting the European options in financial practice, since it reflects the "uncertainty" of the financial market better than the option prices. It is also known that if the stock price drops, the IV raises (so-called leverage effect), motivates hedging strategies based on IVs. Consequently, for the purpose of this application, we will regard the BS–IV as an individual financial variable. The practical relevance of such an approach is justified by the volatility based financial products such as VDAX, which are commonly traded on the option markets.

The goal of this analysis is to study the dynamics of the IV functions for different maturities. More specifically, our aim is to construct low dimensional factor models based on the truncated Karhunen–Loève expansions (1) for the log-returns of the IV functions of options with different maturities and compare these factor models using the methodology presented in the previous sections. Analysis of IVs based on a low-dimensional factor model gives directly a descriptive insight into the structure of distribution of the log-IV-returns—structure of the factors and empirical distribution of the factor loadings may be a good starting point for further pricing models. In practice, such a factor model can also be used in Monte Carlo based pricing methods and for risk-management (hedging) purposes. For comprehensive monographs on IV and IV-factor models, see Hafner (2004) or Fengler (2005b).

The idea of constructing and analyzing the factor models for log-IV-returns for different maturities was originally proposed by Fengler, Härdle and Villa (2003), who studied the dynamics of the IV via PCA on discretized IV functions for different maturity groups and tested the Common Principal Components (CPC) hypotheses (equality of eigenvectors and eigenspaces for different groups). Fengler, Härdle and Villa (2003) proposed a PCA-based factor model for log-IV-returns on (short) maturities 1, 2 and 3 months and grid of moneyness $[0.85, 0.9, 0.95, 1, 1.05, 1.1]$. They showed that the factor functions do not significantly differ and only the factor loadings differ across maturity groups. Their method relies on the CPC methodology introduced by Flury (1988) which is based on maximum likelihood estimation under the assumption of multivariate normality. The log-IV-returns are extracted by the two-dimensional Nadaraya–Watson estimate.

The main aim of this application is to reconsider their results in a functional sense. Doing so, we overcome two basic weaknesses of their approach. First, the factor model proposed by Fengler, Härdle and Villa (2003) is performed only on a sparse design of moneyness. However, in practice (e.g., in Monte Carlo pricing methods), evaluation of the model on a fine grid is needed. Using the functional PCA approach, we may overcome this difficulty and evaluate the factor model on an arbitrary fine grid. The second difficulty of the procedure proposed by Fengler, Härdle and Villa (2003) stems from the data design—on the exchange we cannot observe options with desired maturity on each day and we need to estimate them from the IV-functions with maturities observed on the particular day. Consequently, the two-dimensional Nadaraya–Watson estimator proposed by Fengler, Härdle and Villa (2003) results essentially in the (weighted) average of the IVs (with closest maturities) observed on a particular day, which may affect the test of the common eigenfunction hypothesis. We use the linear interpolation scheme in the *total variance* $\sigma^2_{\mathrm{TOT},i}(\kappa, \tau) \stackrel{\mathrm{def}}{=} (\sigma^\tau_i(\kappa))^2 \tau$, in order to recover the IV functions with fixed maturity (on day $i$). This interpolation scheme is based on the arbitrage arguments originally proposed by Kahalé (2004) for zero-dividend and zero-interest rate case and generalized for deterministic interest rate by Fengler (2005a). More precisely, having IVs with maturities observed on a particular day $i$: $\tilde{\sigma}_i^{\tau_{j_i}}(\kappa)$, $j_i = 1, \ldots, p_{\tau_i}$, we calculate the corresponding total variance $\tilde{\sigma}_{\mathrm{TOT},i}(\kappa, \tau_{j_i})$. From these total variances we linearly interpolate the total variance with the desired maturity from the nearest maturities observed on day $i$. The total variance can be easily transformed to corresponding IV $\tilde{\sigma}_i^\tau(\kappa)$. As the last step, we calculate the log-returns $\Delta \log \tilde{\sigma}_i^\tau(\kappa) \overset{\text{def}}{=} \log \tilde{\sigma}_{i+1}^\tau(\kappa) - \log \tilde{\sigma}_i^\tau(\kappa)$. The log-IV-returns are observed for each maturity $\tau$ on a discrete grid $\kappa_{ik}^\tau$. We assume that observed log-IV-return $\Delta \log \tilde{\sigma}_i^\tau(\kappa_{ik}^\tau)$ consists of true log-return of the IV function denoted by $\Delta \log \sigma_i^\tau(\kappa_{ik}^\tau)$ and possibly of some additional error $\varepsilon_{ik}^\tau$. By setting $Y_{ik}^\tau := \Delta \log \tilde{\sigma}_i^\tau(\kappa_{ik}^\tau)$, $X_i^\tau(\kappa) := \Delta \log \sigma_i^\tau(\kappa)$, we obtain an analogue of the model (4) with the argument $\kappa$:

$$(15) \qquad Y_{ik}^\tau = X_i^\tau(\kappa_{ik}) + \varepsilon_{ik}^\tau, \qquad i = 1, \ldots, n_\tau.$$

In order to simplify the notation and make the connection with the theoretical part clear, we will use the notation of (15).

For our analysis we use a recent data set containing daily data from January 2004 to June 2004 from the German–Swiss exchange (EUREX). Violations of the arbitrage-free assumptions ("obvious" errors in data) were corrected using the procedure proposed by Fengler (2005a). Similarly to Fengler, Härdle and Villa (2003), we excluded options with maturity smaller then 10 days, since these option-prices are known to be very noisy, partially because of a special and arbitrary setup in the pricing systems of the dealers. Using the interpolation scheme described above, we calculate the log-IV-returns for two maturity groups: "1M" group with maturity $\tau = 0.12$ (measured in years) and "3M" group with maturity $\tau = 0.36$. The observed log-IV-returns are denoted by $Y_{ik}^{1M}$, $k = 1, \ldots, K_i^{1M}$, $Y_{ik}^{3M}$, $k = 1, \ldots, K_i^{3M}$. Since we ensured that for no $i$, the interpolation procedure uses data with the same maturity for both groups, this procedure has no impact on the independence of both samples.

The underlying models based on the truncated version of (3) are as follows:

$$(16) \qquad X_i^{1M}(\kappa) = \bar{X}^{1M}(\kappa) + \sum_{r=1}^{L_{1M}} \hat{\beta}_{ri}^{1M} \widehat{\gamma}_r^{1M}(\kappa), \qquad i = 1, \ldots, n_{1M},$$

$$(17) \qquad X_i^{3M}(\kappa) = \bar{X}^{3M}(\kappa) + \sum_{r=1}^{L_{3M}} \hat{\beta}_{ri}^{3M} \widehat{\gamma}_r^{3M}(\kappa), \qquad i = 1, \ldots, n_{3M}.$$

Models (16) and (17) can serve, for example, in a Monte Carlo pricing tool in the risk management for pricing exotic options where the whole path of implied volatilities is needed to determine the price. Estimating the factor functions in (16) and (17) by eigenfunctions displayed in Figure 1, we only need to fit the (estimated) factor loadings $\hat{\beta}_{ji}^{1M}$ and $\hat{\beta}_{ji}^{3M}$. The pillar of the model is the dimension reduction. Keeping the factor function fixed for a certain time period, we need to analyze (two) multivariate random processes of the factor loadings. For the purposes of this paper we will focus on the comparison of factors from models (16) and (17) and the technical details of the factor loading analysis will not be discussed here, since in this respect we refer to Fengler, Härdle and Villa (2003), who proposed to fit the factor loadings by centered normal distributions with diagonal variance matrix containing the corresponding eigenvalues. For a deeper discussion of the fitting of factor loadings using a more sophisticated approach, basically based on (possibly multivariate) GARCH models; see Fengler (2005b).

From our data set we obtained 88 functional observations for the 1M group ($n_{1M}$) and 125 observations for the 3M group ($n_{3M}$). We will estimate the model on the interval for futures moneyness $\kappa \in [0.8, 1.1]$. In comparison to Fengler, Härdle and Villa (2003), we may estimate models (16) and (17) on an arbitrary fine grid (we used an equidistant grid of 500 points on the interval $[0.8, 1.1]$). For illustration, the Nadaraya–Watson (NW) estimator of resulting log-returns is plotted in Figure 2. The smoothing parameters have been chosen in accordance with the requirements in Section 2.2. As argued in Section 2.2, we should use small smoothing parameters in order to avoid a possible bias in the estimated eigenfunctions. Thus, we use for each $i$ essentially the smallest bandwidth $b_i$ that guarantees that estimator $\hat{X}_i$ is defined on the entire support $[0.8, 1.1]$.

Using the procedures described in Section 2.1, we first estimate the eigenfunctions of both maturity groups. The estimated eigenfunctions are plotted in Figure 1. The structure of the eigenfunctions is in accordance with other empirical studies on IV-surfaces. For a deeper discussion and economical interpretation, see, for example, Fengler, Härdle and Mammen (2007) or Fengler, Härdle and Villa (2003).

Clearly, the ratio of the variance explained by the $k$th factor function is given by the quantity $\hat{\nu}_k^{1M} = \hat{\lambda}_k^{1M} / \sum_{j=1}^{n_{1M}} \hat{\lambda}_j^{1M}$ for the 1M group and, correspondingly, by $\hat{\nu}_k^{3M}$ for the 3M group. In Table 3 we list the contributions of the factor functions. Looking at Table 3, we can see that 4th factor functions explain less than 1% of the variation. This number was the "threshold" for the choice of $L_{1M}$ and $L_{2M}$.

We can observe (see Figure 1) that the factor functions for both groups are similar. Thus, in the next step we use the bootstrap test for testing the equality of the factor functions. We use 2000 bootstrap replications. The test of equality of the eigenfunctions was rejected for the first eigenfunction for the analyzed time period (January 2004–June 2004) at a significance level $\alpha = 0.05$ (P-value 0.01). We may conclude that the (first) factor functions are not identical in the factor model for both maturity groups. However, from a practical point of view, we are more interested in checking the appropriateness of the entire models for a fixed number of factors: $L = 2$ or $L = 3$ in (16) and (17). This requirement translates into the testing of the equality of eigenspaces. Thus, in the next step we use the same setup (2000 bootstrap replications) to test the hypotheses that the first two and first three eigenfunctions span the same eigenspaces $\mathcal{E}_L^{1M}$ and $\mathcal{E}_L^{3M}$. None of the hypotheses for $L = 2$ and $L = 3$ is rejected at significance level $\alpha = 0.05$ (P-value is 0.61 for $L = 2$ and 0.09 for $L = 3$). Summarizing, even in the functional sense we have no significant reason to reject the hypothesis of common eigenspaces for these two maturity groups. Using this hypothesis, the factors governing the movement of the returns of IV surface are invariant to time to maturity, only their relative importance can vary. This leads to the common factor model: $X_i^\tau(\kappa) = \bar{X}^\tau(\kappa) + \sum_{r=1}^{L_\tau} \hat{\beta}_{ri}^\tau \widehat{\gamma_r}(\kappa), i = 1, \ldots, n_\tau, \ \tau = 1M, 3M$, where $\gamma_r := \gamma_r^{1M} = \gamma_r^{3M}$. Beside contributing to the understanding of the structure of the IV function dynamics, the common factor model helps us to reduce the number of functional factors by half compared to models (16) and (17). Furthermore, from the technical point of view, we also obtain an additional dimension reduction and higher estimation precision, since under this hypothesis we may estimate the eigenfunctions from the (individually centered) pooled sample $X_i(\kappa)^{1M}, i = 1, \ldots, n_{1M}, X_i^{3M}(\kappa), i =$ $1, \ldots, n_{3M}$. The main improvement compared to the multivariate study by Fengler, Härdle and Villa (2003) is that our test is performed in the functional sense – it does not depend on particular discretization and our factor model can be evaluated on an arbitrary fine grid.

TABLE 3
*Variance explained by the eigenfunctions*

|  | Var. explained 1M | Var. explained 3M |
|---|---|---|
| $\hat{\nu}_1^\tau$ | 89.9% | 93.0% |
| $\hat{\nu}_2^\tau$ | 7.7% | 4.2% |
| $\hat{\nu}_3^\tau$ | 1.7% | 1.0% |
| $\hat{\nu}_4^\tau$ | 0.6% | 0.4% |

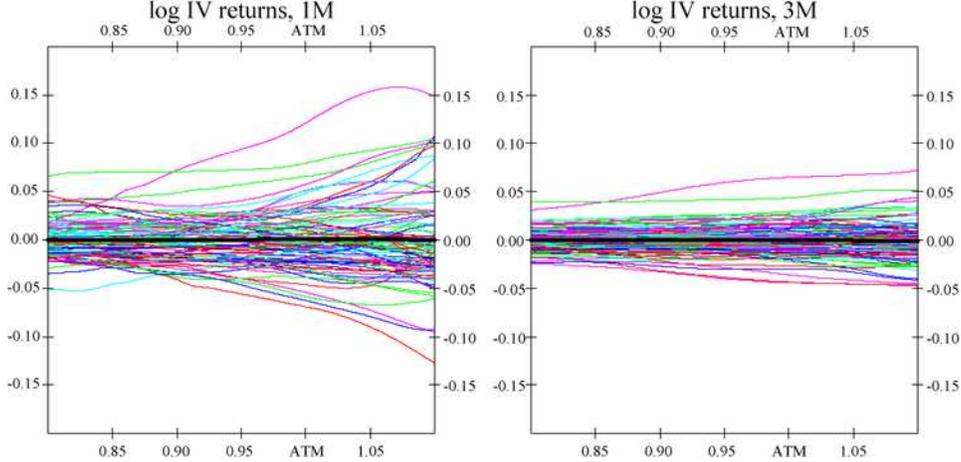

FIG. 2. *Nadaraya–Watson estimate of the log-IV-returns for maturity 1M (left figure) and 3M (right figure). The bold line is the sample mean of the corresponding group.*

## APPENDIX: MATHEMATICAL PROOFS

In the following, $\|v\| = (\int_0^1 v(t)^2\, dt)^{1/2}$ will denote the $L^2$-norm for any square integrable function $v$. At the same time, $\|a\| = (\frac{1}{k}\sum_{i=1}^k a_i^2)^{1/2}$ will indicate the Euclidean norm, whenever $a \in \mathbb{R}^k$ is a $k$-vector for some $k \in \mathbb{N}$.

In the proof of Theorem 1, $\mathrm{E}_\varepsilon$ and $\mathrm{Var}_\varepsilon$ denote expectation and variance with respect to $\varepsilon$ only (i.e., conditional on $t_{ij}$ and $X_i$).

PROOF OF THEOREM 1.  Recall the definition of the $\chi_i(t)$ and note that $\chi_i(t) = \chi_i^X(t) + \chi_i^\varepsilon(t)$, where

$$\chi_i^\varepsilon(t) = \sum_{j=1}^{T_i} \varepsilon_{i(j)} I\left(t \in \left[\frac{t_{i(j-1)} + t_{i(j)}}{2}, \frac{t_{i(j)} + t_{i(j+1)}}{2}\right)\right),$$

as well as

$$\chi_i^X(t) = \sum_{j=1}^{T_i} X_i(t_{i(j)}) I\left(t \in \left[\frac{t_{i(j-1)} + t_{i(j)}}{2}, \frac{t_{i(j)} + t_{i(j+1)}}{2}\right)\right)$$

for $t \in [0,1]$, $t_{i(0)} = -t_{i(1)}$ and $t_{i(T_i+1)} = 2 - t_{i(T_i)}$. Similarly, $\chi_i^*(t) = \chi_i^{X*}(t) + \chi_i^{\varepsilon*}(t)$.

By Assumption 2, $\mathrm{E}(|t_{i(j)} - t_{i(j-1)}|^s) = \mathcal{O}(T^{-s})$ for $s = 1, \ldots, 4$, and the convergence is uniform in $j < n$. Our assumptions on the structure of $X_i$ together with some straightforward Taylor expansions then lead to

$$\langle \chi_i, \chi_j \rangle = \langle X_i, X_j \rangle + \mathcal{O}_p(1/T)$$

and

$$\langle \chi_i, \chi_i^* \rangle = \|X_i\|^2 + \mathcal{O}_p(1/T).$$

Moreover,

$$\mathrm{E}_\varepsilon(\langle \chi_i^\varepsilon, \chi_j^X \rangle) = 0, \qquad\qquad \mathrm{E}_\varepsilon(\|\chi_i^\varepsilon\|^2) = \sigma_i^2,$$

$$\mathrm{E}_\varepsilon(\langle \chi_i^\varepsilon, \chi_i^{*} \rangle) = 0, \qquad\qquad \mathrm{E}_\varepsilon(\langle \chi_i^\varepsilon, \chi_i^{\varepsilon*} \rangle^2) = \mathcal{O}_p(1/T),$$

$$\mathrm{E}_\varepsilon(\langle \chi_i^\varepsilon, \chi_j^X \rangle^2) = \mathcal{O}_p(1/T), \qquad \mathrm{E}_\varepsilon(\langle \chi_i^\varepsilon, \chi_j^X \rangle \langle \chi_k^\varepsilon, \chi_l^X \rangle) = 0 \qquad \text{for } i \neq k,$$

$\mathrm{E}_\varepsilon(\langle \chi_i^\varepsilon, \chi_j^\varepsilon \rangle \langle \chi_i^\varepsilon, \chi_k^\varepsilon \rangle) = 0$ \quad for $j \neq k$ and $\mathrm{E}_\varepsilon(\|\chi_i^\varepsilon\|^4) = \mathcal{O}_p(1)$

hold (uniformly) for all $i, j = 1, \ldots, n$.

Consequently, $\mathrm{E}_\varepsilon(\|\bar\chi\|^2 - \|\bar X\|^2) = \mathcal{O}_p(T^{-1} + n^{-1})$.

When using these relations, it is easily seen that for all $i, j = 1, \ldots, n$

$$
\widehat{M}_{ij} - M_{ij} = \mathcal{O}_p(T^{-1/2} + n^{-1}) \quad \text{and}
$$

(18)

$$
\mathrm{tr}\{(\widehat{M} - M)^2\}^{1/2} = \mathcal{O}_p(1 + nT^{-1/2}).
$$

Since the orthonormal eigenvectors $p_q$ of $M$ satisfy $\|p_q\| = 1$, we furthermore obtain for any $i = 1, \ldots, n$ and all $q = 1, 2, \ldots$

$$
(19) \qquad \sum_{j=1}^{n} p_{jq} \left\{ \widehat{M}_{ij} - M_{ij} - \int_0^1 \chi_i^\varepsilon(t) \chi_j^X(t)\, dt \right\} = \mathcal{O}_p(T^{-1/2} + n^{-1/2}),
$$

as well as

$$
(20) \qquad \sum_{j=1}^{n} p_{jq} \int_0^1 \chi_i^\varepsilon(t) \chi_j^X(t)\, dt = \mathcal{O}_p\left( \frac{n^{1/2}}{T^{1/2}} \right)
$$

and

$$
(21) \qquad \sum_{i=1}^{n} a_i \sum_{j=1}^{n} p_{jq} \int_0^1 \chi_i^\varepsilon(t) \chi_j^X(t)\, dt = \mathcal{O}_p\left( \frac{n^{1/2}}{T^{1/2}} \right)
$$

for any further vector $a$ with $\|a\| = 1$.

Recall that the $j$th largest eigenvalue $l_j$ satisfies $n\hat{\lambda}_j = l_j$. Since by assumption $\inf_{s \neq r} |\lambda_r - \lambda_s| > 0$, the results of Dauxois, Pousse and Romain (1982) imply that $\hat{\lambda}_r$ converges to $\lambda_r$ as $n \to \infty$, and $\sup_{s \neq r} \frac{1}{|\lambda_r - \lambda_s|} = \mathcal{O}_p(1)$, which leads to $\sup_{s \neq r} \frac{1}{|l_r - l_s|} = \mathcal{O}_p(1/n)$. Assertion (a) of Lemma A of Kneip and Utikal (2001) together with (18)–(21) then implies that

$$
\left| \hat{\lambda}_r - \frac{\hat{l}_r}{n} \right| = n^{-1} |l_r - \hat{l}_r| = n^{-1} |p_r^\top (\widehat{M} - M) p_r| + \mathcal{O}_p(T^{-1} + n^{-1})
$$

(22)

$$
= \mathcal{O}_p\{(nT)^{-1/2} + T^{-1} + n^{-1}\}.
$$

When analyzing the difference between the estimated and true eigenvectors $\hat{p}_r$ and $p_r$, assertion (b) of Lemma A of Kneip and Utikal (2001) together with (18) lead to

$$
(23) \qquad \hat{p}_r - p_r = -\mathcal{S}_r(\widehat{M} - M) p_r + \mathcal{R}_r, \qquad \text{with } \|\mathcal{R}_r\| = \mathcal{O}_p(T^{-1} + n^{-1})
$$

and $\mathcal{S}_r = \sum_{s \neq r} \frac{1}{l_s - l_r} p_s p_s^\top$. Since $\sup_{\|a\|=1} a^\top \mathcal{S}_r a \leq \sup_{s \neq r} \frac{1}{|l_r - l_s|} = \mathcal{O}_p(1/n)$, we can conclude that

$$
(24) \qquad \|\hat{p}_r - p_r\| = \mathcal{O}_p(T^{-1/2} + n^{-1}),
$$

and our assertion on the sequence $n^{-1} \sum_i (\hat{\beta}_{ri} - \hat{\beta}_{ri;T})^2$ is an immediate consequence.

Let us now consider assertion (ii). The well-known properties of local linear estimators imply that $|\mathrm{E}_\varepsilon\{\hat{X}_i(t) - X_i(t)\}| = \mathcal{O}_p(b^2)$, as well as $\mathrm{Var}_\varepsilon\{\hat{X}_i(t)\} = \mathcal{O}_p\{Tb\}$, and the convergence is uniform for all $i, n$. Furthermore, due to the independence of the error term $\varepsilon_{ij}$, $\mathrm{Cov}_\varepsilon\{\hat{X}_i(t), \hat{X}_j(t)\} = 0$ for $i \neq j$. Therefore,

$$\left|\hat{\gamma}_r(t) - \frac{1}{\sqrt{l_r}}\sum_{i=1}^n p_{ir}\hat{X}_i(t)\right| = \mathcal{O}_p\left(b^2 + \frac{1}{\sqrt{nTb}}\right).$$

On the other hand, (18)–(24) imply that with $\hat{X}(t) = (\hat{X}_1(t), \ldots, \hat{X}_n(t))^\top$

$$
\begin{aligned}
&\left|\hat{\gamma}_{r;T}(t) - \frac{1}{\sqrt{l_r}}\sum_{i=1}^n p_{ir}\hat{X}_i(t)\right| \\
&\quad = \left|\frac{1}{\sqrt{l_r}}\sum_{i=1}^n(\hat{p}_{ir} - p_{ir})X_i(t) + \frac{1}{\sqrt{l_r}}\sum_{i=1}^n(\hat{p}_{ir} - p_{ir})\{\hat{X}_i(t) - X_i(t)\}\right| \\
&\qquad + \mathcal{O}_p(T^{-1} + n^{-1}) \\
&\quad = \frac{\|\mathcal{S}_r X(t)\|}{\sqrt{l_r}}\left|p_r^\top(\hat{M} - M)\mathcal{S}_r\frac{X(t)}{\|\mathcal{S}_r X(t)\|}\right| \\
&\qquad + \mathcal{O}_p(b^2 T^{-1/2} + T^{-1}b^{-1/2} + n^{-1}) \\
&\quad = \mathcal{O}_p(n^{-1/2}T^{-1/2} + b^2 T^{-1/2} + T^{-1}b^{-1/2} + n^{-1}).
\end{aligned}
$$

This proves the theorem. $\square$

PROOF OF THEOREM 2. First consider assertion (i). By definition,

$$\bar{X}(t) - \mu(t) = n^{-1}\sum_{i=1}^n\{X_i(t) - \mu(t)\} = \sum_r\left(n^{-1}\sum_{i=1}^n\beta_{ri}\right)\gamma_r(t).$$

Recall that, by assumption, $\beta_{ri}$ are independent, zero mean random variables with variance $\lambda_r$, and that the above series converges with probability 1. When defining the truncated series

$$V(q) = \sum_{r=1}^q\left(n^{-1}\sum_{i=1}^n\beta_{ri}\right)\gamma_r(t),$$

standard central limit theorems therefore imply that $\sqrt{n}V(q)$ is asymptotically $N(0, \sum_{r=1}^q\lambda_r\gamma_r(t)^2)$ distributed for any possible $q \in \mathbb{N}$.

The assertion of a $N(0, \sum_{r=1}^\infty\lambda_r\gamma_r(t)^2)$ limiting distribution now is a consequence of the fact that for all $\delta_1, \delta_2 > 0$ there exists a $q_\delta$ such that $P\{|\sqrt{n}V(q) - \sqrt{n}\sum_r(n^{-1}\sum_{i=1}^n\beta_{ri})\gamma_r(t)| > \delta_1\} < \delta_2$ for all $q \geq q_\delta$ and all $n$ sufficiently large.

In order to prove assertions (i) and (ii), consider some fixed $r \in \{1, 2, \ldots\}$ with $\lambda_{r-1} > \lambda_r > \lambda_{r+1}$. Note that $\Gamma$ as well as $\hat{\Gamma}_n$ are nuclear, self-adjoint and non-negative linear operators with $\Gamma v = \int \sigma(t, s) v(s) \, ds$ and $\hat{\Gamma}_n v = \int \hat{\sigma}(t, s) v(s) \, ds$, $v \in L^2[0, 1]$. For $m \in \mathbb{N}$, let $\Pi_m$ denote the orthogonal projector from $L^2[0, 1]$ into the $m$-dimensional linear space spanned by $\{\gamma_1, \ldots, \gamma_m\}$, that is, $\Pi_m v = \sum_{j=1}^{m} \langle v, \gamma_j \rangle \gamma_j$, $v \in L^2[0, 1]$. Now consider the operator $\Pi_m \hat{\Gamma}_n \Pi_m$, as well as its eigenvalues and corresponding eigenfunctions denoted by $\hat{\lambda}_{1,m} \geq \hat{\lambda}_{2,m} \geq \cdots$ and $\hat{\gamma}_{1,m}, \hat{\gamma}_{2,m}, \ldots$, respectively. It follows from well-known results in the Hilbert space theory that $\Pi_m \hat{\Gamma}_n \Pi_m$ converges strongly to $\hat{\Gamma}_n$ as $m \to \infty$. Furthermore, we obtain (Rayleigh–Ritz theorem)

$$(25) \quad \lim_{m \to \infty} \hat{\lambda}_{r,m} = \lambda_r \quad \text{and} \quad \lim_{m \to \infty} \|\hat{\gamma}_r - \hat{\gamma}_{r,m}\| = 0 \quad \text{if } \hat{\lambda}_{r-1} > \hat{\lambda}_r > \hat{\lambda}_{r+1}.$$

Note that under the above condition $\hat{\gamma}_r$ is uniquely determined up to sign, and recall that we always assume that the right "versions" (with respect to sign) are used so that $\langle \hat{\gamma}_r, \hat{\gamma}_{r,m} \rangle \geq 0$. By definition, $\beta_{ji} = \int \gamma_j(t) \{X_i(t) - \mu(t)\} \, dt$, and therefore, $\int \gamma_j(t) \{X_i(t) - \bar{X}(t)\} \, dt = \beta_{ji} - \bar{\beta}_j$, as well as $X_i - \bar{X} = \sum_j (\beta_{ji} - \bar{\beta}_j) \gamma_j$, where $\bar{\beta}_j = \frac{1}{n} \sum_{i=1}^{n} \beta_{ji}$. When analyzing the structure of $\Pi_m \hat{\Gamma}_n \Pi_m$ more deeply, we can verify that $\Pi_m \hat{\Gamma}_n \Pi_m v = \int \hat{\sigma}_m(t, s) v(s) \, ds$, $v \in L^2[0, 1]$, with

$$\hat{\sigma}_m(t, s) = g_m(t)^\top \hat{\Sigma}_m g_m(s),$$

where $g_m(t) = (\gamma_1(t), \ldots, \gamma_m(t))^\top$, and where $\hat{\Sigma}_m$ is the $m \times m$ matrix with elements $\{\frac{1}{n} \sum_{i=1}^{n} (\beta_{ji} - \bar{\beta}_j)(\beta_{ki} - \bar{\beta}_k)\}_{j,k=1,\ldots,m}$. Let $\lambda_1(\hat{\Sigma}_m) \geq \lambda_2(\hat{\Sigma}_m) \geq \cdots \geq \lambda_m(\hat{\Sigma}_m)$ and $\hat{\zeta}_{1,m}, \ldots, \hat{\zeta}_{m,m}$ denote eigenvalues and corresponding eigenvectors of $\hat{\Sigma}_m$. Some straightforward algebra then shows that

$$(26) \quad \hat{\lambda}_{r,m} = \lambda_r(\hat{\Sigma}_m), \qquad \hat{\gamma}_{r,m} = g_m(t)^\top \hat{\zeta}_{r,m}.$$

We will use $\Sigma_m$ to represent the $m \times m$ diagonal matrix with diagonal entries $\lambda_1 \geq \cdots \geq \lambda_m$. Obviously, the corresponding eigenvectors are given by the $m$-dimensional unit vectors denoted by $e_{1,m}, \ldots, e_{m,m}$. Lemma A of [Kneip and Utikal (2001)](#) now implies that the differences between eigenvalues and eigenvectors of $\Sigma_m$ and $\hat{\Sigma}_m$ can be bounded by

$$(27) \quad \hat{\lambda}_{r,m} - \lambda_r = \mathrm{tr}\{e_{r,m} e_{r,m}^\top (\hat{\Sigma}_m - \Sigma_m)\} + \tilde{R}_{r,m},$$

$$\text{with } \tilde{R}_{r,m} \leq \frac{6 \sup_{\|a\|=1} a^\top (\hat{\Sigma}_m - \Sigma_m)^2 a}{\min_s |\lambda_s - \lambda_r|},$$

$$(28) \quad \hat{\zeta}_{r,m} - e_{r,m} = -S_{r,m}(\hat{\Sigma}_m - \Sigma_m) e_{r,m} + R_{r,m}^*,$$

$$\text{with } \|R_{r,m}^*\| \leq \frac{6 \sup_{\|a\|=1} a^\top (\hat{\Sigma}_m - \Sigma_m)^2 a}{\min_s |\lambda_s - \lambda_r|^2},$$

where $S_{r,m} = \sum_{s \neq r} \frac{1}{\lambda_s - \lambda_r} e_{s,m} e_{s,m}^\top$.

Assumption 1 implies $\mathrm{E}(\bar{\beta}_r) = 0$, $\mathrm{Var}(\bar{\beta}_r) = \frac{\lambda_r}{n}$, and with $\delta_{ii} = 1$, as well as $\delta_{ij} = 0$ for $i \neq j$, we obtain

$$
\begin{aligned}
\mathrm{E}\bigg\{ &\sup_{\|a\|=1} a^\top (\hat{\Sigma}_m - \Sigma_m)^2 a \bigg\} \\
&\leq \mathrm{E}\{\mathrm{tr}[(\hat{\Sigma}_m - \Sigma_m)^2]\} \\
&= \mathrm{E}\bigg\{ \sum_{j,k=1}^m \bigg[ \frac{1}{n} \sum_{i=1}^n (\beta_{ji} - \bar{\beta}_j)(\beta_{ki} - \bar{\beta}_k) - \delta_{jk} \lambda_j \bigg]^2 \bigg\} \\
&\leq \mathrm{E}\bigg\{ \sum_{j,k=1}^\infty \bigg[ \frac{1}{n} \sum_{i=1}^n (\beta_{ji} - \bar{\beta}_j)(\beta_{ki} - \bar{\beta}_k) - \delta_{jk} \lambda_j \bigg]^2 \bigg\} \\
&= \frac{1}{n} \bigg( \sum_j \sum_k \mathrm{E}\{\beta_{ji}^2 \beta_{ki}^2\} \bigg) + \mathcal{O}(n^{-1}) = \mathcal{O}(n^{-1}),
\end{aligned}
$$

(29)

for all $m$. Since $\mathrm{tr}\{e_{r,m} e_{r,m}^\top (\hat{\Sigma}_m - \Sigma_m)\} = \frac{1}{n} \sum_{i=1}^n (\beta_{ri} - \bar{\beta}_r)^2 - \lambda_r$, (25), (26), (27) and (29) together with standard central limit theorems imply that

$$
\begin{aligned}
\sqrt{n}(\hat{\lambda}_r - \lambda_r) &= \frac{1}{\sqrt{n}} \sum_{i=1}^n (\beta_{ri} - \bar{\beta}_r)^2 - \lambda_r + \mathcal{O}_p(n^{-1/2}) \\
&= \frac{1}{\sqrt{n}} \sum_{i=1}^n [(\beta_{ri})^2 - \mathrm{E}\{(\beta_{ri})^2\}] + \mathcal{O}_p(n^{-1/2}) \\
&\xrightarrow{\mathcal{L}} N(0, \Lambda_r).
\end{aligned}
$$

(30)

It remains to prove assertion (iii). Relations (26) and (28) lead to

$$
\begin{aligned}
\hat{\gamma}_{r,m}(t) - \gamma_r(t) &= g_m(t)^\top (\hat{\zeta}_{r,m} - e_{r,m}) \\
&= -\sum_{s \neq r} \bigg\{ \frac{1}{n(\lambda_s - \lambda_r)} \sum_{i=1}^n (\beta_{si} - \bar{\beta}_s)(\beta_{ri} - \bar{\beta}_r) \bigg\} \gamma_s(t) \\
&\quad + g_m(t)^\top R_{r,m}^*,
\end{aligned}
$$

(31)

where due to (29) the function $g_m(t)^\top R_{r,m}^*$ satisfies

$$
\begin{aligned}
\mathrm{E}(\|g_m^\top R_{r,m}^*\|) &= \mathrm{E}(\|R_{r,m}^*\|) \\
&\leq \frac{6}{n \min_s |\lambda_s - \lambda_r|^2} \bigg( \sum_j \sum_k \mathrm{E}\{\beta_{ji}^2 \beta_{ki}^2\} \bigg) + \mathcal{O}(n^{-1}),
\end{aligned}
$$

for all $m$. By Assumption 1, the series in (31) converge with probability 1 as $m \to \infty$.

Obviously, the event $\hat{\lambda}_{r-1} > \hat{\lambda}_r > \hat{\lambda}_{r+1}$ occurs with probability 1. Since $m$ is arbitrary, we can therefore conclude from (25) and (31) that

$$
\begin{aligned}
(32) \quad & \hat{\gamma}_r(t) - \gamma_r(t) \\
& = -\sum_{s \neq r} \left\{ \frac{1}{n(\lambda_s - \lambda_r)} \sum_{i=1}^{n} (\beta_{si} - \bar{\beta}_s)(\beta_{ri} - \bar{\beta}_r) \right\} \gamma_s(t) \; + \; R_r^*(t) \\
& = -\sum_{s \neq r} \left\{ \frac{1}{n(\lambda_s - \lambda_r)} \sum_{i=1}^{n} \beta_{si} \beta_{ri} \right\} \gamma_s(t) \; + \; R_r(t),
\end{aligned}
$$

where $\|R_r^*\| = \mathcal{O}_p(n^{-1})$, as well as $\|R_r\| = \mathcal{O}_p(n^{-1})$. Moreover, $\sqrt{n} \times \sum_{s \neq r} \{ \frac{1}{n(\lambda_s - \lambda_r)} \sum_{i=1}^{n} \beta_{si} \beta_{ri} \} \gamma_s(t)$ is a zero mean random variable with variance $\sum_{q \neq r} \sum_{s \neq r} \frac{\mathrm{E}[\beta_{ri}^2 \beta_{qi} \beta_{si}]}{(\lambda_q - \lambda_r)(\lambda_s - \lambda_r)} \gamma_q(t) \gamma_s(t) < \infty$. By Assumption 1, it follows from standard central limit arguments that for any $q \in \mathbb{N}$ the truncated series $\sqrt{n} W(q) \stackrel{\text{def}}{=} \sqrt{n} \sum_{s=1, s \neq r}^{q} [\frac{1}{n(\lambda_s - \lambda_r)} \sum_{i=1}^{n} \beta_{si} \beta_{ri}] \gamma_s(t)$ is asymptotically normal distributed. The asserted asymptotic normality of the complete series then follows from an argument similar to the one used in the proof of assertion (i).  □

Proof of Theorem 3.  The results of Theorem 2 imply that

$$
\begin{aligned}
(33) \quad n\Delta_1 = \int \Bigg( & \sum_r \frac{1}{\sqrt{q_1 n_1}} \sum_{i=1}^{n_1} \beta_{ri}^{(1)} \gamma_r^{(1)}(t) \\
& - \sum_r \frac{1}{\sqrt{q_2 n_2}} \sum_{i=1}^{n_2} \beta_{ri}^{(2)} \gamma_r^{(2)}(t) \Bigg)^2 dt.
\end{aligned}
$$

Furthermore, independence of $X_i^{(1)}$ and $X_i^{(2)}$ together with (30) imply that

$$
\begin{aligned}
(34) \quad & \sqrt{n}[\hat{\lambda}_r^{(1)} - \lambda_r^{(1)} - \{\hat{\lambda}_r^{(2)} - \lambda_r^{(2)}\}] \stackrel{\mathcal{L}}{\to} N\left(0, \frac{\Lambda_r^{(1)}}{q_1} + \frac{\Lambda_r^{(2)}}{q_2}\right) \quad \text{and} \\
& \frac{n}{\Lambda_r^{(1)}/q_1 + \Lambda_r^{(2)}/q_2} \Delta_{3,r} \stackrel{\mathcal{L}}{\to} \chi_1^2.
\end{aligned}
$$

Furthermore, (32) leads to

$$
\begin{aligned}
(35) \quad n\Delta_{2,r} = \Bigg\| & \sum_{s \neq r} \left\{ \frac{1}{\sqrt{q_1 n_1}(\lambda_s^{(1)} - \lambda_r^{(1)})} \sum_{i=1}^{n_1} \beta_{si}^{(1)} \beta_{ri}^{(1)} \right\} \gamma_s^{(1)} \\
& - \sum_{s \neq r} \left\{ \frac{1}{\sqrt{q_2 n_2}(\lambda_s^{(2)} - \lambda_r^{(2)})} \sum_{i=1}^{n_2} \beta_{si}^{(2)} \beta_{ri}^{(2)} \right\} \gamma_s^{(2)} \Bigg\|^2 + \mathcal{O}_p(n^{-1/2})
\end{aligned}
$$

and

$$n\Delta_{4,L} = n \iint \Bigg[ \sum_{r=1}^{L} \gamma_r^{(1)}(t)\{\hat{\gamma}_r^{(1)}(u) - \gamma_r^{(1)}(u)\}$$

$$+ \gamma_r^{(1)}(u)\{\hat{\gamma}_r^{(1)}(t) - \gamma_r^{(1)}(t)\}$$

$$- \sum_{r=1}^{L} \gamma_r^{(2)}(t)\{\hat{\gamma}_r^{(2)}(u) - \gamma_r^{(2)}(u)\}$$

$$+ \gamma_r^{(2)}(u)\{\hat{\gamma}_r^{(2)}(t) - \gamma_r^{(2)}(t)\} \Bigg]^2 dt\,du + \mathcal{O}_p(n^{-1/2})$$

$$(36) \qquad = \iint \Bigg[ \sum_{r=1}^{L} \sum_{s>L} \Bigg\{ \frac{1}{\sqrt{q_1 n_1}(\lambda_s^{(1)} - \lambda_r^{(1)})} \sum_{i=1}^{n_1} \beta_{si}^{(1)}\beta_{ri}^{(1)} \Bigg\}$$

$$\times \{\gamma_r^{(1)}(t)\gamma_s^{(1)}(u) + \gamma_r^{(1)}(u)\gamma_s^{(1)}(t)\}$$

$$- \sum_{r=1}^{L} \sum_{s>L} \Bigg\{ \frac{1}{\sqrt{q_2 n_2}(\lambda_s^{(2)} - \lambda_r^{(2)})} \sum_{i=1}^{n_2} \beta_{si}^{(2)}\beta_{ri}^{(2)} \Bigg\}$$

$$\times \{\gamma_r^{(2)}(t)\gamma_s^{(2)}(u) + \gamma_r^{(2)}(u)\gamma_s^{(2)}(t)\} \Bigg]^2 dt\,du$$

$$+ \mathcal{O}_p(n^{-1/2}).$$

In order to verify (36), note that $\sum_{r=1}^{L} \sum_{s=1, s\neq r}^{L} \frac{1}{(\lambda_s^{(p)} - \lambda_r^{(p)})} a_r a_s = 0$ for $p = 1, 2$ and all possible sequences $a_1, \ldots, a_L$. It is clear from our assumptions that all sums involved converge with probability 1. Recall that $\mathrm{E}(\beta_{ri}^{(p)}\beta_{si}^{(p)}) = 0$, $p = 1, 2$ for $r \neq s$.

It follows that $\tilde{X}_r^{(p)} := \frac{1}{\sqrt{q_p n_p}} \sum_{s\neq r} \sum_{i=1}^{n_p} \frac{\beta_{si}^{(p)}\beta_{ri}^{(p)}}{\lambda_s^{(p)} - \lambda_r^{(p)}} \gamma_s^{(p)}$, $p = 1, 2$, is a continuous, zero mean random function on $L^2[0,1]$, and, by assumption, $\mathrm{E}(\|\tilde{X}_r^{(p)}\|^2) < \infty$. By Hilbert space central limit theorems [see, e.g., [Araujo and Giné (1980)](#)], $\tilde{X}_r^{(p)}$ thus converges in distribution to a Gaussian random function $\xi_r^{(p)}$ as $n \to \infty$. Obviously, $\xi_r^{(1)}$ is independent of $\xi_r^{(2)}$. We can conclude that $n\Delta_{4,L}$ possesses a continuous limit distribution $F_{4,L}$ defined by the distribution of $\iint [\sum_{r=1}^{L} \{\xi_r^{(1)}(t)\gamma_r^{(1)}(u) + \xi_r^{(1)}(u)\gamma_r^{(1)}(t)\} - \sum_{r=1}^{L} \{\xi_r^{(2)}(t)\gamma_r^{(2)}(u) + \xi_r^{(2)}(u) \times \gamma_r^{(2)}(t)\}]^2 dt\,du$. Similar arguments show the existence of continuous limit distributions $F_1$ and $F_{2,r}$ of $n\Delta_1$ and $n\Delta_{2,r}$.

For given $q \in \mathbb{N}$, define vectors $b_{i1}^{(p)} = (\beta_{1i}^{(p)}, \ldots, \beta_{qi}^{(p)},)^\top \in \mathbb{R}^q$, $b_{i2}^{(p)} = (\beta_{1i}^{(p)}\beta_{ri}^{(p)}, \ldots, \beta_{r-1,i}^{(p)}\beta_{ri}^{(p)}, \beta_{r+1,i}^{(p)}\beta_{ri}^{(p)}, \ldots, \beta_{qi}^{(p)}\beta_{ri}^{(p)})^\top \in \mathbb{R}^{q-1}$ and $b_{i3} = (\beta_{1i}^{(p)}\beta_{2i}^{(p)},$ $\dots, \beta_{qi}^{(p)} \beta_{Li}^{(p)})^\top \in \mathbb{R}^{(q-1)L}$. When the infinite sums over $r$ in (33), respectively $s \neq r$ in (35) and (36), are restricted to $q \in \mathbb{N}$ components (i.e., $\sum_r$ and $\sum_{s>L}$ are replaced by $\sum_{r \leq q}$ and $\sum_{L<s\leq q}$), then the above relations can generally be presented as limits $n\Delta = \lim_{q\to\infty} n\Delta(q)$ of quadratic forms

$$n\Delta_1(q) = \begin{pmatrix} \dfrac{1}{\sqrt{n_1}}\displaystyle\sum_{i=1}^{n_1} b_{i1}^{(1)} \\ \dfrac{1}{\sqrt{n_2}}\displaystyle\sum_{i=1}^{n_2} b_{i1}^{(2)} \end{pmatrix}^\top Q_1^q \begin{pmatrix} \dfrac{1}{\sqrt{n_1}}\displaystyle\sum_{i=1}^{n_1} b_{i1}^{(1)} \\ \dfrac{1}{\sqrt{n_2}}\displaystyle\sum_{i=1}^{n_2} b_{i1}^{(2)} \end{pmatrix},$$

$$(37) \qquad n\Delta_{2,r}(q) = \begin{pmatrix} \dfrac{1}{\sqrt{n_1}}\displaystyle\sum_{i=1}^{n_1} b_{i2}^{(1)} \\ \dfrac{1}{\sqrt{n_2}}\displaystyle\sum_{i=1}^{n_2} b_{i2}^{(2)} \end{pmatrix}^\top Q_2^q \begin{pmatrix} \dfrac{1}{\sqrt{n_1}}\displaystyle\sum_{i=1}^{n_1} b_{i2}^{(1)} \\ \dfrac{1}{\sqrt{n_2}}\displaystyle\sum_{i=1}^{n_2} b_{i2}^{(2)} \end{pmatrix},$$

$$n\Delta_{4,L}(q) = \begin{pmatrix} \dfrac{1}{\sqrt{n_1}}\displaystyle\sum_{i=1}^{n_1} b_{i3}^{(1)} \\ \dfrac{1}{\sqrt{n_2}}\displaystyle\sum_{i=1}^{n_2} b_{i3}^{(2)} \end{pmatrix}^\top Q_3^q \begin{pmatrix} \dfrac{1}{\sqrt{n_1}}\displaystyle\sum_{i=1}^{n_1} b_{i3}^{(1)} \\ \dfrac{1}{\sqrt{n_2}}\displaystyle\sum_{i=1}^{n_2} b_{i3}^{(2)} \end{pmatrix},$$

where the elements of the $2q \times 2q$, $2(q-1) \times 2(q-1)$ and $2L(q-1) \times 2L(q-1)$ matrices $Q_1^q$, $Q_2^q$ and $Q_3^q$ can be computed from the respective ($q$-element) version of (33)–(36). Assumption 1 implies that all series converge with probability 1 as $q \to \infty$, and by (33)–(36), it is easily seen that for all $\epsilon, \delta > 0$ there exist some $q(\epsilon,\delta), n(\epsilon,\delta) \in \mathbb{N}$ such that

$$(38) \qquad P(|n\Delta_1 - n\Delta_1(q)| > \epsilon) < \delta, \qquad P(|n\Delta_{2,r} - n\Delta_{2,r}(q)| > \epsilon) < \delta,$$
$$P(|n\Delta_{4,L} - n\Delta_{4,L}(q)| > \epsilon) < \delta$$

hold for all $q \geq q(\epsilon,\delta)$ and all $n \geq n(\epsilon,\delta)$. For any given $q$, we have $\mathrm{E}(b_{i1}) = \mathrm{E}(b_{i2}) = \mathrm{E}(b_{i3}) = 0$, and it follows from Assumption 1 that the respective covariance structures can be represented by finite covariance matrices $\Omega_{1,q}$, $\Omega_{2,q}$ and $\Omega_{3,q}$. It therefore follows from our assumptions together with standard multivariate central limit theorems that the vectors $\{\frac{1}{\sqrt{n_1}}\sum_{i=1}^{n_1}(b_{ik}^{(1)})^\top,$ $\frac{1}{\sqrt{n_2}}\sum_{i=1}^{n_2}(b_{ik}^{(2)})^\top\}^\top$, $k = 1, 2, 3$, are asymptotically normal with zero means and covariance matrices $\Omega_{1,q}$, $\Omega_{2,q}$ and $\Omega_{3,q}$. One can thus conclude that, as $n \to \infty$,

$$(39) \qquad n\Delta_1(q) \xrightarrow{\mathcal{L}} F_{1,q}, \qquad n\Delta_{2,r}(q) \xrightarrow{\mathcal{L}} F_{2,r,q}, \qquad n\Delta_{4,L}(q) \xrightarrow{\mathcal{L}} F_{4,L,q},$$

where $F_{1,q}, F_{2,r,q}, F_{4,L,q}$ denote the continuous distributions of the quadratic forms $z_1^\top Q_1^q z_1$, $z_2^\top Q_2^q z_2$, $z_3^\top Q_3^q z_3$ with $z_1 \sim N(0, \Omega_{1,q})$, $z_2 \sim N(0, \Omega_{2,q})$, $z_3 \sim$

$N(0, \Omega_{3,q})$. Since $\epsilon, \delta$ are arbitrary, (38) implies

$$(40) \qquad \lim_{q \to \infty} F_{1,q} = F_1, \qquad \lim_{q \to \infty} F_{2,r,q} = F_{2,r}, \qquad \lim_{q \to \infty} F_{4,L,q} = F_{4,L}.$$

We now have to consider the asymptotic properties of bootstrapped eigenvalues and eigenfunctions. Let $\bar{X}^{(p)*} = \frac{1}{n_p} \sum_{i=1}^{n_p} X_i^{(p)*}$, $\beta_{ri}^{(p)*} = \int \gamma_r^{(p)}(t) \{X_i^{(p)*}(t) - \mu(t)\}$, $\bar{\beta}_r^{(p)*} = \frac{1}{n_p} \sum_{i=1}^{n_p} \beta_{ri}^{(p)*}$, and note that $\int \gamma_r^{(p)}(t) \{X_i^{(p)*}(t) - \bar{X}^{(p)*}(t)\} = \beta_{ri}^{(p)*} - \bar{\beta}_r^{(p)*}$. When considering unconditional expectations, our assumptions imply that for $p = 1, 2$

$$\mathrm{E}[\beta_{ri}^{(p)*}] = 0, \qquad\qquad \mathrm{E}[(\beta_{ri}^{(p)*})^2] = \lambda_r^{(p)},$$

$$\mathrm{E}[(\bar{\beta}_r^{(p)*})^2] = \frac{\lambda_r^{(p)}}{n_p}, \qquad \mathrm{E}\{[(\beta_{ri}^{(p)*})^2 - \lambda_r^{(p)}]^2\} = \Lambda_r^{(p)},$$

$$(41) \qquad \mathrm{E}\left\{ \sum_{l,k=1}^{\infty} \left[ \frac{1}{n_p} \sum_{i=1}^{n_p} (\beta_{li}^{(p)*} - \bar{\beta}_l^{(p)*})(\beta_{ki}^{(p)*} - \bar{\beta}_k^{(p)*}) - \delta_{lk} \lambda_l^{(p)} \right]^2 \right\}$$

$$= \frac{1}{n_p} \left( \sum_l \Lambda_l^{(p)} + \sum_{l \neq k} \lambda_l^{(p)} \lambda_k^{(p)} \right) + \mathcal{O}(n_p^{-1}).$$

One can infer from (41) that the arguments used to prove Theorem 1 can be generalized to approximate the difference between the bootstrap eigenvalues and eigenfunctions $\hat{\lambda}_r^{(p)*}$, $\hat{\gamma}_r^{(p)*}$ and the true eigenvalues $\lambda_r^{(p)}$, $\gamma_r^{(p)}$. All infinite sums involved converge with probability 1. Relation (30) then generalizes to

$$\sqrt{n_p}(\hat{\lambda}_r^{(p)*} - \hat{\lambda}_r^{(p)})$$

$$= \sqrt{n_p}(\hat{\lambda}_r^{(p)*} - \lambda_r^{(p)}) - \sqrt{n_p}(\hat{\lambda}_r^{(p)} - \lambda_r^{(p)})$$

$$(42) \qquad = \frac{1}{\sqrt{n_p}} \sum_{i=1}^{n_p} (\beta_{ri}^{(p)*} - \bar{\beta}_r^{(p)*})^2$$

$$\qquad\qquad - \frac{1}{\sqrt{n_p}} \sum_{i=1}^{n_p} (\beta_{ri}^{(p)} - \bar{\beta}_r^{(p)})^2 + \mathcal{O}_p(n_p^{-1/2})$$

$$= \frac{1}{\sqrt{n_p}} \sum_{i=1}^{n_p} \left\{ (\beta_{ri}^{(p)*})^2 - \frac{1}{n_p} \sum_{k=1}^{n_p} (\beta_{rk}^{(p)})^2 \right\} + \mathcal{O}_p(n_p^{-1/2}).$$

Similarly, (32) becomes

$$\hat{\gamma}_r^{(p)*} - \hat{\gamma}_r^{(p)}$$

$$(43) \qquad = \hat{\gamma}_r^{(p)*} - \gamma_r^{(p)} - (\hat{\gamma}_r^{(p)} - \gamma_r^{(p)})$$

$$= -\sum_{s \neq r} \Bigg\{ \frac{1}{\lambda_s^{(p)} - \lambda_r^{(p)}} \frac{1}{n_p} \sum_{i=1}^{n_p} (\beta_{si}^{(p)*} - \bar{\beta}_s^{(p)*})(\beta_{ri}^{(p)*} - \bar{\beta}_r^{(p)*})$$

$$- \frac{1}{\lambda_s^{(p)} - \lambda_r^{(p)}} \frac{1}{n_p} \sum_{i=1}^{n_p} (\beta_{si}^{(p)} - \bar{\beta}_s^{(p)})(\beta_{ri}^{(p)} - \bar{\beta}_r^{(p)}) \Bigg\} \gamma_s^{(p)}(t)$$

$$+ R_r^{(p)*}(t)$$

$$= -\sum_{s \neq r} \Bigg\{ \frac{1}{\lambda_s^{(p)} - \lambda_r^{(p)}} \frac{1}{n_p} \sum_{i=1}^{n_p} \left( \beta_{si}^{(p)*} \beta_{ri}^{(p)*} - \frac{1}{n_p} \sum_{k=1}^{n_p} \beta_{sk}^{(p)} \beta_{rk}^{(p)} \right) \Bigg\} \gamma_s^{(p)}(t)$$

$$+ \tilde{R}_r^{(p)*}(t),$$

where due to (28), (29) and (41), the remainder term satisfies $\|R_r^{(p)*}\| = \mathcal{O}_p(n_p^{-1})$.

We are now ready to analyze the bootstrap versions $\Delta^*$ of the different $\Delta$. First consider $\Delta_{3,r}^*$ and note that $\{(\beta_{ri}^{(p)*})^2\}$ are i.i.d. bootstrap resamples from $\{(\beta_{ri}^{(p)})^2\}$. It therefore follows from basic bootstrap results that the conditional distribution of $\frac{1}{\sqrt{n_p}} \sum_{i=1}^{n_p} [(\beta_{ri}^{(p)*})^2 - \frac{1}{n_p} \sum_{k=1}^{n_p} (\beta_{rk}^{(p)})^2]$ given $\mathcal{X}_p$ converges to the same $N(0, \Lambda_r^{(p)})$ limit distribution as $\frac{1}{\sqrt{n_p}} \sum_{i=1}^{n_p} [(\beta_{ri}^{(p)})^2 - \mathrm{E}\{(\beta_{ri}^{(p)})^2\}]$. Together with the independence of $(\beta_{ri}^{(1)*})^2$ and $(\beta_{ri}^{(2)*})^2$, the assertion of the theorem is an immediate consequence.

Let us turn to $\Delta_1^*$, $\Delta_{2,r}^*$ and $\Delta_{4,L}^*$. Using (41)–(43), it is then easily seen that $n\Delta_1^*$, $n\Delta_{2,r}^*$ and $n\Delta_{4,L}^*$ admit expansions similar to (33), (35) and (36), when replacing there $\frac{1}{\sqrt{n_p}} \sum_{i=1}^{n_p} \beta_{ri}^{(p)}$ by $\frac{1}{\sqrt{n_p}} \sum_{i=1}^{n_p} (\beta_{ri}^{(p)*} - \frac{1}{n_p} \sum_{k=1}^{n_p} \beta_{rk}^{(p)})$, as well as $\frac{1}{\sqrt{n_p}} \sum_{i=1}^{n_p} \beta_{si}^{(p)} \beta_{ri}^{(p)}$ by $\frac{1}{\sqrt{n_p}} \sum_{i=1}^{n_p} (\beta_{si}^{(p)*} \beta_{ri}^{(p)*} - \frac{1}{n_p} \sum_{k=1}^{n_p} \beta_{sk}^{(p)} \beta_{rk}^{(p)})$.

Replacing $\beta_{ri}^{(p)}$, $\beta_{si}^{(p)}$ by $\beta_{ri}^{(p)*}$, $\beta_{si}^{(p)*}$ leads to bootstrap analogs $b_{ik}^{(p)*}$ of the vectors $b_{ik}^{(p)}$, $k = 1, 2, 3$. For any $q \in \mathbb{N}$, define bootstrap versions $n\Delta_1^*(q)$, $n\Delta_{2,r}^*(q)$ and $n\Delta_{4,L}^*(q)$ of $n\Delta_1(q)$, $n\Delta_{2,r}(q)$ and $n\Delta_{4,L}(q)$ by using $(\frac{1}{\sqrt{n_1}} \sum_{i=1}^{n_1} (b_{ik}^{(1)*} - \frac{1}{n_1} \sum_{i=1}^{n_1} b_{ik}^{(1)})^\top, \frac{1}{\sqrt{n_2}} \sum_{i=1}^{n_2} (b_{ik}^{(2)*} - \frac{1}{n_2} \sum_{k=1}^{n_2} b_{ik}^{(2)})^\top)$ instead of $(\frac{1}{\sqrt{n_1}} \sum_{i=1}^{n_1} (b_{ik}^{(1)})^\top, \frac{1}{\sqrt{n_2}} \sum_{i=1}^{n_2} (b_{ik}^{(2)})^\top)$, $k = 1, 2, 3$, in (37). Applying again (41)–(43), one can conclude that for any $\epsilon > 0$ there exists some $q(\epsilon)$ such that, as $n \to \infty$,

$$(44) \qquad \begin{aligned} P(|n\Delta_1^* - n\Delta_1^*(q)| < \epsilon) &\to 1, \\ P(|n\Delta_{2,r}^* - n\Delta_{2,r}^*(q)| < \epsilon) &\to 1, \\ P(|n\Delta_{4,L}^* - n\Delta_{4,L}^*(q)| < \epsilon) &\to 1 \end{aligned}$$

hold for all $q \geq q(\epsilon)$. Of course, (44) generalizes to the conditional probabilities given $\mathcal{X}_1$, $\mathcal{X}_2$.

In order to prove the theorem, it thus only remains to show that for *any* given $q$ and all $\delta$

$$(45) \qquad |\mathrm{P}(n\Delta(q) \geq \delta) - \mathrm{P}(n\Delta^*(q) \geq \delta \mid \mathcal{X}_1, \mathcal{X}_2)| = o_p(1)$$

hold for either $\Delta(q) = \Delta_1(q)$ and $\Delta^*(q) = \Delta_1^*(q)$, $\Delta(q) = \Delta_{2,r}(q)$ and $\Delta^*(q) = \Delta_{2,r}^*(q)$, or $\Delta(q) = \Delta_{4,L}(q)$ and $\Delta^*(q) = \Delta_{4,L}^*(q)$. But note that for $k = 1, 2, 3$, $\mathrm{E}(b_{ik}) = 0$, $\{b_{ik}^{(j)*}\}$ are i.i.d. bootstrap resamples from $\{b_{ik}^{(p)}\}$, and $\mathrm{E}(b_{ik}^{(p)*} \mid \mathcal{X}_1, \mathcal{X}_2) = \frac{1}{n_p} \sum_{k=1}^{n_p} b_{ik}^{(p)}$ are the corresponding conditional means. It therefore follows from basic bootstrap results that as $n \to \infty$ the conditional distribution of $(\frac{1}{\sqrt{n_1}} \sum_{i=1}^{n_1} (b_{ik}^{(1)*} - \frac{1}{n_1} \sum_{k=1}^{n_1} b_{ik}^{(1)})^\top, \frac{1}{\sqrt{n_2}} \sum_{i=1}^{n_2} (b_{ik}^{(2)*} - \frac{1}{n_2} \sum_{k=1}^{n_2} b_{ik}^{(2)})^\top)$ given $\mathcal{X}_1$, $\mathcal{X}_2$ converges to the same $N(0, \Omega_{k,q})$ limit distribution as $(\frac{1}{\sqrt{n_1}} \sum_{i=1}^{n_1} (b_{ik}^{(1)})^\top, \frac{1}{\sqrt{n_2}} \sum_{i=1}^{n_2} (b_{ik}^{(2)})^\top)$. This obviously holds for all $q \in \mathbb{N}$, and (45) is an immediate consequence. The theorem then follows from (38), (39), (40), (44) and (45). $\quad \square$

M. BENKO
W. HÄRDLE
CASE—CENTER FOR APPLIED STATISTICS AND ECONOMICS
HUMBOLDT-UNIVERSITÄT ZU BERLIN
SPANDAUERSTR 1
D-10178 BERLIN
GERMANY
E-MAIL: benko@wiwi.hu-berlin.de
        haerdle@wiwi.hu-berlin.de
URL: http://www.case.hu-berlin.de/

A. KNEIP
STATISTISCHE ABTEILUNG
DEPARTMENT OF ECONOMICS
UNIVERSITÄT BONN
ADENAUERALLEE 24-26
D-53113 BONN
GERMANY
E-MAIL: akneip@uni-bonn.de